\documentclass[final,3p]{elsarticle}

\usepackage{amsmath}
\usepackage{amsthm}
\usepackage{amssymb}
\usepackage{graphicx}
\usepackage{mathrsfs}
\usepackage{color}
\usepackage{mathtools} % for mathclap
\usepackage{booktabs}% for top mid and bottomrule in tables
\usepackage{epsfig} % for epsilon figures
\usepackage{subcaption} % for subcaptionbox environment
\usepackage{float} % For [H] to force position of figure
\usepackage{multirow}
\usepackage[bookmarks,bookmarksnumbered]{hyperref}
\hypersetup{colorlinks = true,linkcolor = blue,anchorcolor =red,citecolor = blue,filecolor = red,urlcolor = red,
            pdfauthor=author}
            
\newtheorem{remark}{Remark}

\DeclareMathAlphabet\mathbfcal{OMS}{cmsy}{b}{n} % Bold calligraphic mathbfcal
\newcommand{\bs}{\mathbf}

\begin{document}

\begin{frontmatter}


%% and the form \Ex\ead[url] for the home page:
%% \Ex\ead{email address}
%% \Ex\ead[url]{home page}

\title{Data-free Non-intrusive Model Reduction for Nonlinear Finite Element Models via Spectral Submanifolds}

\author[label1]{Mingwu Li}
\author[label2]{Thomas Thurnher}
\author[label2]{Zhenwei Xu}
\author[label3]{Shobhit Jain\corref{cor1}}
\ead{Shobhit.Jain@tudelft.nl}
\cortext[cor1]{Corresponding author}
\address[label1]{Department of Mechanics and Aerospace Engineering, Southern University of Science and Technology, Shenzhen, China}
\address[label2]{Institute for Mechanical Systems, ETH Z\"{u}rich Leonhardstrasse 21, 8092 Z\"{u}rich, Switzerland}
\address[label3]{Delft Institute of Applied Mathematics, TU Delft, Mekelweg 4, 2628CD, Delft, The Netherlands}

\begin{abstract}
The theory of spectral submanifolds (SSMs) has emerged as a powerful tool for constructing rigorous, low-dimensional reduced-order models (ROMs) of high-dimensional nonlinear mechanical systems. A direct computation of SSMs requires explicit knowledge of nonlinear coefficients in the equations of motion, which limits their applicability to generic finite-element (FE) solvers. Here, we propose a non-intrusive algorithm for the computation of the SSMs and the associated ROMs up to arbitrary polynomial orders. This non-intrusive algorithm only requires system nonlinearity as a black box and hence, enables SSM-based model reduction via generic finite-element software. Our expressions and algorithms are valid for systems with up to cubic-order nonlinearities, including velocity-dependent nonlinear terms, asymmetric damping, and stiffness matrices, and hence work for a large class of mechanics problems. We demonstrate the effectiveness of the proposed non-intrusive approach over a variety of FE examples of increasing complexity, including a micro-resonator FE model containing more than a million degrees of freedom.
\end{abstract}


\end{frontmatter}

\section{Introduction}

% spectral submanifolds

Linear modal subspaces are powerful constructs for rigorous model reduction of linear time-invariant systems. Due to the invariance of modal subspaces, a full simulation trajectory with initialized on a modal subspace stays on it for all times. For instance, consider a weakly damped linear mechanical system, initialized at rest in the shape of a vibration mode. This system will freely vibrate with decaying amplitudes while retaining the initial shape due to the invariance of the modal subspaces associated with the system's vibration modes. In addition, a slow modal subspace associated with the slowest decaying modes (low-frequency modes) will attract nearby solution trajectories exponentially fast. As a result, one can project the equations governing the full system onto a low-dimensional slow modal subspace for rigorous model reduction of linear systems. Indeed, the success of modal truncation in constructing lower-dimensional reduced-order models (ROMs) for high-dimensional linear systems~\cite{hintz1975analytical,dickens1997critique} can be justified based on this principle. Given this success, a natural question is, can we extend the idea of modal truncation to nonlinear problems?

The key ingredients of the above extension are \emph{invariance} and attractivity or \emph{normal hyperbolicity}. Nonlinear normal modes (NNMs) have been defined as invariant manifolds that serve as nonlinear continuations of modal subspaces. Model reduction via NNMs enables the prediction of free and forced response in nonlinear vibrations~\cite{Shaw1993,Pesheck2002AManifolds,jiang2005nonlinear}. There are, however, infinitely many such invariant manifolds for a given modal subspace. The smoothest invariant manifold among these exists uniquely under appropriate non-resonant conditions and is called as a spectral submanifold (SSM)~\cite{ssmexist,haller2023nonlinear}. In particular, a slow SSM, which is the extension of the slow modal subspace, is also an attracting invariant manifold and hence can be used to perform exact model reduction for nonlinear problems~\cite{ssmexist,haller2023nonlinear}.

In the recent years, SSMs have emerged as accurate and predictive model reduction tools for high-dimensional dynamical systems~\cite{Ponsioen2020,Jain2021HowModels,part-i,Cenedese2022Data-drivenSubmanifolds}. For systems without internal resonances, reduction on two-dimensional SSMs enables analytical extraction of backbone curves and forced response curves and surfaces~\cite{Breunung2018ExplicitSystems,li2024fast}. For harmonically excited systems with internal resonances, one can solve for fixed points of the SSM-based ROMs to predict nonlinear periodic response of the full system~\cite{part-i}. Likewise, a limit cycle of the SSM-based ROMs corresponds to a torus on which a quasi-periodic orbit of the full system stays~\cite{part-i}. Further, bifurcation analysis of the forced periodic and quasi-periodic orbits can be simplified as the bifurcation detection of the fixed point and limit cycle of SSM-based ROMs~\cite{part-ii,li2023nonlinear}. More recently, SSM-based model reduction has also been extended to systems with parametric excitation~\cite{thurnher2023nonautonomous} and to constrained mechanical systems, where the equations of motion are in the form of differential-algebraic equations~\cite{li2023model}.

% SSM computational methods
SSMs and their associated ROMs can be computed in an automated fashion using the parameterization method~\cite{CabreI,CabreIII,Haro2016TheManifolds}. In earlier studies~\cite{ponsioen2018automated,Breunung2018ExplicitSystems,Ponsioen2020}, such computations are based on equations of motion in modal basis, where the coefficient matrix is a diagonal matrix, with each diagonal entry being an eigenvalue of the linear part of the dynamics. These computations are ineffective for high-dimensional problems because one must first transform the equations of motion into modal coordinates. This transformation needs all eigenvectors (modes) of the linear part of the dynamics. Indeed, as illustrated in~\cite{Jain2021HowModels}, this transformation will destroy the sparsity of nonlinear terms, and also, computing all eigenvectors becomes prohibitive as the dimension increases. In recent developments~\cite{Jain2021HowModels,part-i,Vizzaccaro2021HighOD,thurnher2023nonautonomous}, SSM computations were directly conducted in physical coordinates using only the eigenvectors associated with the master subspace of the SSM. This improvement makes it possible to compute SSMs of high-dimensional FE models.

% why do we need non-intrusive computation
The above computations of SSMs, however, are fall in the category of intrusive reduction because they need explicit access to the nonlinear terms in the equations of motion. These terms  are typically inaccessible via generic or commercial FE software. \emph{Non-intrusive} or \emph{indirect} reduction methods, on the other hand, use the generic FE solvers with minimal access to their internal core and have a compelling use-case as they are agnostic to the choice of the solver~\cite{Mignolet2013AStructures}. For projection-based model reduction, the Stiffness Evaluation Procedure (STEP)~\cite{Muravyov2003DeterminationStructures} and Implicit Condensation and Expansion (ICE))~\cite{mcewan2001finite,Hollkamp2008Reduced-orderExpansion} are prominent nonintrusive techniques that extract the nonlinear coefficients for equations of motion in modal coordinates. In the STEP method, a set of prescribed displacements, which are linear combinations of modal basis, are carefully designed and then imposed on the structure to extract the nonlinear coefficients~\cite{Muravyov2003DeterminationStructures}. In the ICE method, a set of prescribed external loads that are also related to the linear combinations of modal basis is applied to the structure, and the resulting static displacements are solved and then projected to modal basis, followed by a linear regression to fit the nonlinear coefficients~\cite{Hollkamp2008Reduced-orderExpansion}.

%The STEP and ICE methods need to be further combined with model reduction schemes because of the curse of dimensionality for high-dimensional FE problems~\cite{Mignolet2013AStructures,Touze2021ModelTechniques}. Indeed, computing all eigenvectors and extracting all nonlinear coefficients are infeasible for high-dimensional FE models. STEP combined with modal truncation was used in earlier studies to construct ROMs~\cite{Muravyov2003DeterminationStructures}. 

The convergence of such modal ROMs with increasing number of modes, however, is problem-dependent and can be slow in practice~\cite{Givois2019OnModels,Vizzaccaro2020Non-intrusiveElements}. To speed up this convergence, dual modes~\cite{wang2009nonlinear,kim2013nonlinear} and modal derivatives~\cite{mahdiabadi2021non} may be included in the STEP framework for projection-based model reduction. Still, projection-based reduction methods are generally not justifiable for nonlinear systems due to the lack of invariance of linear subspaces. Consequently, their performance is guaranteed only when a slow-fast decomposition is satisfied~\cite{haller2017exact,Shen2021ReducedApproach,Vizzaccaro2021ComparisonDerivatives}. Such a decomposition, however, generally does not hold for complex engineering structures~\cite{part-i}. A natural question then is if we can combine non-intrusive procedures such as STEP~\cite{Muravyov2003DeterminationStructures} with the computation of attracting invariant-manifolds for rigorous nonlinear model reduction. 

% In ICE method, the membrane coordinates are parameterized via bending coordinates, and then membrane basis functions are constructed, resulting in reduction by static condensation~\cite{Shen2021ReducedApproach}

% contributions of this study 

An important contribution~\cite{Vizzaccaro2021DirectStructures} proposed formal computation of direct normal forms up to cubic order for non-intrusive model reduction of FE problems. This method~\cite{Vizzaccaro2021DirectStructures} provides explicit expressions for the normal form over a selection of modes using the STEP~\cite{Muravyov2003DeterminationStructures}, albeit only up to cubic order. While higher-order expansions have been considered intrusively in later developments~\cite{Vizzaccaro2021HighOD,opreni2023high}, their non-intrusive implementations remain unexplored. Furthermore, these contributions~\cite{Vizzaccaro2021DirectStructures,Vizzaccaro2021HighOD,opreni2023high} develop expressions for proportionally damped systems and with purely geometric nonlinearities. Indeed, these damping assumptions simplify nonintrusive implementations by evaluating the nonlinearity over purely real inputs. Generally damped mechanical systems, however, exhibit complex eigenvectors, which makes nonintrusive implementations via the STEP challenging due to the evaluation of the nonlinearity over complex inputs. This is especially the case for commerical FE software, where nonlinearity evaluations are optimized for real inputs and complex inputs may result in incorrect function evaluations. The above limitations motivate our current study.

% Reduction via the normal form method and STEP has been applied to geometrically nonlinear flat structures~\cite{Givois2019OnModels}. However, as demonstrated in~\cite{Givois2019OnModels,Vizzaccaro2020Non-intrusiveElements}, the coupling coefficients characterizing modal interactions converge very slowly as the number of modal basis of FE model increases, resulting in intensive computations of STEP. Indeed, the slow convergence becomes more obvious for 3D FE problems~\cite{Vizzaccaro2020Non-intrusiveElements}. To address the slow convergence regarding the number of mode basis, 

\subsection*{Our contributions} 
We develop a non-intrusive procedure for computing spectral submanifolds and their reduced dynamics up to arbitrary polynomial orders in nonlinear FE models. Our expressions are valid for high-dimensional mechanical systems with up to cubic-order nonlinearities including velocity-dependent nonlinear forces, e.g., in viscoelasticity with large deformations~\cite{amabili2018nonlinear,amabili2019derivation}, and with possibly asymmetric damping and stiffness matrices, which are common in rotating machinery and fluid-structure interaction with gyroscopic and follower forces~\cite{xie2020structural,mereles2023model,paidoussis1998fluid,li2023model}. Furthermore, we show how our nonintrusive routine can be coupled with generic FE software, including commercial packages that are not meant to accept complex-valued inputs, e.g., COMSOL Multiphysics~\cite{multiphysics1998introduction}.  We implement our nonintrusive routine in the open-source package, SSMTool~\cite{SSMTool2} and demonstrate its applicability and effectiveness towards nonlinear forced response predictions over a variety of FE applications. Our demonstrations include examples such as a shell structure with internal resonance, an aircraft wing, a viscoelastic cover plate, fluid structure interaction in a pipe conveying fluid, and a 3D continuum-based FE model of a MEMS gyroscope.

Thanks to the STEP~\cite{Muravyov2003DeterminationStructures}, our non-intrusive model reduction via SSMs is \emph{data-free}, in contrast to the existing \emph{data-assisted} or \emph{data-driven} procedures used for nontrusive SSM-based model reduction~\cite{Kaszas2022,cenedese2023data}. Data-based non-intrusive computations require full system simulations to generate unforced trajectory that is used to fit nonlinear parts of the SSM and its reduced dynamics parametrizations. For high-dimensional FE models, the computational cost of these full system simulations can be prohibitively expensive. Our approach addresses this limitation of data-based techniques as it does not involve any full-system simulations.

The rest of this paper is organized below. We present a review of SSM theory and the intrusive computational methods for SSM-based reduction in Sect.~\ref{sec:ssm-review}. We then show how to compute the SSMs and their associated ROMs in a non-intrusive fashion in Sect.~\ref{sec:SemiIntrusive}. Subsequently, we provide the implementation details of the proposed non-intrusive algorithm in Sect.~\ref{sec:implementation}. In Sect.~\ref{sec:examples}, we demonstrate the effectiveness of the non-intrusive routine with a suite of examples with increasing complexity. We conclude this study in Sect.~\ref{sec:conclusion}.

\section{Review of SSM theory and computation}
\label{sec:ssm-review}
\subsection{Setup and SSM theory}
Consider a harmonically excited nonlinear mechanical system; the equations of motion are 
\begin{equation}
\label{eq:2nd-ode}
\mathbf{M\ddot{x} + C}\mathbf{\dot{x} + Kx} + \mathbf{f}(\mathbf{x},\mathbf{\dot{x}}) = \epsilon \mathbf{f}^\mathrm{ext}({\Omega} t),\quad 0\leq\epsilon\ll1,
\end{equation}
where $\bs{x}\in\mathbb{R}^n$ is a vector of generalized displacement, $\dot{\bs{x}}$ and $\ddot{\bs{x}}$ are vectors of generalized velocity and acceleration, $\bs{M}$, $\bs{C}$, and $\bs{K}$ are mass, damping, and stiffness matrices, $\bs{f}$ represent nonlinear internal forces upto cubic orders, and $\bs{f}^\mathrm{ext}$ denotes external harmonic excitation. We allow for asymmetric mass, damping, and stiffness matrices, and also for velocity-dependent nonlinearities. We assume that the system nonlinearities are upto cubic order, which is sufficient for a wide range of FE-based engineering applications. Let $\bs{z}=(\bs{x},\dot{\bs{x}})$, the equations of motion~\eqref{eq:2nd-ode} can be rewritten as a first-order form below
\begin{equation}\label{eq:DS_naut}
    \mathbf{B}\dot{\mathbf{z}} = \mathbf{Az} + \mathbf{F}(\mathbf{z})  +\epsilon\mathbf{F}^\mathrm{ext}(\mathbf{\Omega}t,\mathbf{z}),
\end{equation}
where
\begin{gather} 
    \mathbf{A} =  \begin{bmatrix} -\mathbf{K} & \mathbf{0} \\ \mathbf{0} & \mathbf{M} \end{bmatrix}, \quad
    \mathbf{B} =  \begin{bmatrix} \mathbf{C} & \mathbf{M} \\ \mathbf{M} & \mathbf{0} \end{bmatrix},\quad 
    \mathbf{F}(\bs{z}) =  \begin{bmatrix}- \mathbf{f}(\bs{x},\dot{\bs{x}}) \\ \mathbf{0} \end{bmatrix}, \quad
    \mathbf{F}^\mathrm{ext}({\Omega} t) =  \begin{bmatrix} \mathbf{f}^\mathrm{ext}({\Omega} t) \\ \mathbf{0}\end{bmatrix}.\label{eq:ABF}
\end{gather}

Consider the linear part of the dynamical system~\eqref{eq:DS_naut}
\begin{equation}
\label{eq:LinDynSys}
    \mathbf{B}\dot{\mathbf{z}} = \mathbf{Az},
\end{equation}
which leads to a generalized eigenvalue problem below
\begin{equation}\label{eq:Eigenprob_Order1_Right}
    (\mathbf{A} - \lambda_j \mathbf{B}) \mathbf{v}_j = \mathbf{0},\quad  \mathbf{w}^*_j(\mathbf{A} - \lambda_j \mathbf{B})= \mathbf{0}
\end{equation}
for $j=1,\cdots,N$. Here $\lambda_j$ is a generalized eigenvalue, and $\mathbf{v}_j$ and $\mathbf{w}_j$ are the corresponding right and left eigenvectors. We assume that the real parts of all eigenvalues are nonzero such that the fixed point $\bs{z}=\bs{0}$ is asymptotically stable. As demonstrated in~\cite{Jain2021HowModels}, it is computationally expensive to compute the full spectrum of~\eqref{eq:LinDynSys} when $N$ is large. In contrast, small subsets of eigenvalues and their corresponding generalized eigenvectors can be computed efficiently even for high-dimensional systems~\cite{Jain2021HowModels,Golub1996MatrixComputations}.

We consider an even subset of $M\ll N$ eigenvalues and arrange them according to their real parts' increasing order of magnitudes. The corresponding eigenvectors span an invariant subspace
\begin{equation}
\mathcal{E} = \mathrm{span}\{\mathbf{v}_1,\cdots,\mathbf{v}_{M}\}
\end{equation}
to~\eqref{eq:LinDynSys}. This subspace can be used as the \emph{master spectral subspace} for model reduction of~\eqref{eq:LinDynSys}. We normalize these eigenvectors such that $(\mathbf{w}_j)^*\mathbf{B}\mathbf{v}_i = \delta_{ij}$. Accordingly, the spectrum of the master subspace $\mathcal{E}$ is given by
\begin{equation}
\mathrm{spect}(\mathcal{E}) = \{\lambda_1,\cdots,\lambda_{M}\}.
\end{equation}

Under the addition of the nonlinearity $\mathbf{F}(\mathbf{z})$ to~\eqref{eq:LinDynSys}, the flat invariant subspace $\mathcal{E}$ does not remain invariant. Instead, it is perturbed into invariant manifolds tangent to $\mathcal{E}$ at the origin for the autonomous system
\begin{equation}\label{eq:DS_aut}
    \mathbf{B}\dot{\mathbf{z}} = \mathbf{Az} + \mathbf{F}(\mathbf{z}) .
\end{equation}
There are infinitely many such invariant manifolds for a given master subspace $\mathcal{E}$~\cite{ssmexist}. among these invariant manifolds, the smoothest one uniquely exists under generic non-resonance conditions~\cite{ssmexist}. This unique invariant manifold is defined as the \emph{spectral submanifold} (SSM) associated to $\mathcal{E}$, and we denote it as $\mathcal{W}(\mathcal{E})$. The dynamics on such an SSM is called its \emph{reduced dynamics} and serves as a mathematically rigorous reduced-order model for the full nonlinear system.

Under further addition of the small-amplitude periodic forcing $\epsilon\mathbf{f}^\mathrm{ext}({\Omega}t)$, the hyperbolic fixed point is perturbed as periodic orbit $\gamma_\epsilon$ of the same stability type. Meanwhile, the autonomous SSM $\mathcal{W}(\mathcal{E})$ becomes periodic and hence a \emph{non-autonomous} SSM $\mathcal{W}_{\gamma_\epsilon}(\mathcal{E})$, provided that appropriate non-resonance conditions are satisfied~\cite{ssmexist}. Let $\bs{p}\in\mathbb{C}^M$ be a vector of reduced coordinates; the non-autonomous SSM can be parameterized via the map below
\begin{equation}
\mathbf{W}_\epsilon(\mathbf{p},\bs\phi): \mathcal{U} = U \times \mathbb{T} \mapsto \mathbb{R}^{N}, U \subset \mathbb{C}^M,
\end{equation}
and the reduced dynamics on the SSM $\mathcal{W}_{\gamma_\epsilon}(\mathcal{E})$ can be described as
\begin{equation}
\dot{\mathbf{p}} = \mathbf{R}_\epsilon (\mathbf{p},\phi), \quad \dot{\phi} = {\Omega}.
\end{equation}
Importantly, the SSM parameterization and the associated reduced dynamics satisfy the invariance equation below
\begin{equation}\label{eq:InvEq}
\mathbf{B} ((\partial_\mathbf{p} \mathbf{W}_\epsilon) \mathbf{R}_\epsilon + {\Omega} \partial_{\phi} \mathbf{W}_\epsilon)
= \mathbf{A}\mathbf{W}_\epsilon + \mathbf{F}\circ \mathbf{W}_\epsilon + \epsilon \mathbf{F}^\mathrm{ext}(\phi).
\end{equation}

\subsection{Computation of SSMs}

To compute an SSM, we express its parametrization in terms of Taylor expansion in the coefficient $\epsilon$. This is feasible due to the smooth dependence of the SSM $\mathcal{W}_{\gamma_\epsilon}(\mathcal{E})$ on the coefficient $\epsilon$. Specifically, the Taylor expansion gives
\begin{align}\label{eq:ExpandW_e}
\mathbf{W}_\epsilon(\mathbf{p,\bs\phi)} &= \mathbf{W}\mathbf{(p)} + \epsilon \mathbf{X}\mathbf{(p,\bs\phi)} + 					O(\epsilon^2),
\\
\label{eq:ExpandR_e}
 \mathbf{R}_\epsilon\mathbf{(p,\bs\phi)} 
&=
    \mathbf {R}\mathbf{(p)} + \epsilon \mathbf{S} \mathbf{(p,\bs\phi)} + O(\epsilon^2).
\end{align}
These expansions allow us to separate autonomous terms from the non-autonomous terms resulting from the external excitation $\epsilon\mathbf{F}^\mathrm{ext}({\Omega}t)$ in \eqref{eq:DS_naut}~\cite{Jain2021HowModels}. Next, we compute the autonomous and non-autonomous terms separately.

\subsubsection{Autonomous manifold}

We substitute the expansions \eqref{eq:ExpandW_e} and \eqref{eq:ExpandR_e} into the invariance equation \eqref{eq:InvEq} and collect terms at $O(\epsilon^0)$, resulting in
\begin{equation}\label{eq:InvEq_aut}
\mathbf{B} (\text{D} \mathbf{W}) \mathbf{R} = \mathbf{A}\mathbf{W} + \mathbf{F} \circ \mathbf{W}.
\end{equation}
The autonomous SSM parametrization and reduced dynamics are further expanded using multivariate monomials
\begin{gather}
    \mathbf{W}(\mathbf{p}) =  \sum_{\mathbf{m}\in \mathbb{N}^M} \mathbf{W}_{\mathbf{m}} \mathbf{p}^\mathbf{m},
\quad \mathbf{W}_{\mathbf{m}} \in \mathbb{C}^N,  \label{eq:Wexpansion}\\
    \mathbf{R}(\mathbf{p}) =  \sum_{\mathbf{m}\in \mathbb{N}^M} \mathbf{R}_{\mathbf{m}} \mathbf{p}^\mathbf{m},
\quad \mathbf{R}_{\mathbf{m}} \in \mathbb{C}^M \label{eq:Rexpansion}.
\end{gather}
Here $\mathbf{m} \in \mathbb{N}^M$ is a multi-index vector, which enables monomials to be written in a compact form as $\mathbf{p^m}=p_1^{m_1} \cdots p_M^{m_M}$. This multi-index notation removes the redundancy of tensor-based polynomial representation and optimizes memory requirements for computations~\cite{thurnher2023nonautonomous}.

We substitute the expansions~\eqref{eq:Wexpansion}-\eqref{eq:Rexpansion} into \eqref{eq:InvEq_aut}, and collect terms in the monomial $\mathbf{p}^\mathbf{m}$ for a given multi-index vector $\mathbf{m}$, yielding the autonomous invariance equation for the multi-index vector $\mathbf{m}$~\cite{thurnher2023nonautonomous}
\begin{align}\label{app_eq:Hom_Aut_First}
\underbrace{
     \bigg(\mathbf{A}
 -
   \Lambda_{\mathbf{m}} \mathbf{B}\bigg)
   }_{:= \mathbfcal{L}_{\mathbf{m}}}
    \mathbf{W}_\mathbf{m}
 =
\sum_{j=1}^M \mathbf{B} \mathbf{v}_j R^j_\mathbf{m}
+
\mathcal{C}_\mathbf{m}
- \mathbf{F \circ W}\Bigr|_{\mathbf{m}},
\end{align}
where $\Lambda_\mathbf{m} = \sum_{i=1}^{M}\lambda_i m_i$ and
\begin{equation}
    \mathcal{C}_\mathbf{m} :=  
    \sum_{j=1}^M  \sum_{\substack{\mathbf{u}, \mathbf{k} \in \mathbb{N}^M \\ \mathbf{u} + \mathbf{k} - \mathbf{e}_j = \mathbf{m} \\ 1 < |\mathbf{u}| < m }} \mathbf{B}\mathbf{W}_{\mathbf{u}} u_j R^j_{\mathbf{k}}.
\end{equation}
At the linear order $|\mathbf{m}|=1$, the nonlinear function $\mathbf{F}$ does not contribute, and we are left with the eigenvalue problem \eqref{eq:Eigenprob_Order1_Right}. This can be easily solved by choosing $\mathbf{W}_{ \mathbf{e}_i} = \mathbf{v}_i$ and $R^j_{\mathbf{e}_i} = \lambda_i \delta_{ij}$~\cite{Jain2021HowModels}. At higher orders $|\mathbf{m}|>1$, the expansion coefficients can be computed in an iterative procedure, as detailed in~\cite{thurnher2023nonautonomous}.

We use a normal-form-style parameterization~\cite{haro2016parameterization} to solve the system of linear equations~\eqref{app_eq:Hom_Aut_First}. In case of near-resonances of the form
\begin{align}\label{eq:autResCond}
     \Lambda_ \mathbf{m} \approx \lambda_i, \qquad i \in \{ 1, \cdots, M \},
\end{align}
the coefficient matrix $\mathbfcal{L}_\mathbf{m}$ is nearly singular~\cite{Jain2021HowModels,thurnher2023nonautonomous}. To lift the singularity, the reduced dynamics is chosen such that $ \mathbf{w}_i^*\mathbfcal{L}_{\mathbf{m}}=0$, which results in
\begin{align}
    R^i_{\mathbf{m}} =- \mathbf{w}_i^* ( \mathcal{C}_\mathbf{m},
- \mathbf{F \circ W}\Bigr|_{\mathbf{m}} ).
\end{align}
Here we have used the orthonormalization of eigenvectors, i.e., $(\mathbf{w}_j)^*\mathbf{B}\mathbf{v}_i = \delta_{ij}$. In the case that~\eqref{eq:autResCond} does not hold, the coefficient matrix $\mathbfcal{L}_\mathbf{m}$ is regular and we simply set $R^j_{\mathbf{m}}=0$ to obtain reduced dynamics in the simplest form.

\subsubsection{Non-autonomous manifold}
We take a leading-order approximation to the non-autonomous part of the SSM, namely,
\begin{equation}
    \bs{X}(\bs{p},\phi)\approx\bs{X}_\bs{0}(\phi),\quad \bs{S}(\bs{p},\phi)\approx\bs{S}_\bs{0}(\phi).
\end{equation}
Substituting the leading order approximation above into the invariance equation~\eqref{eq:InvEq} at $\mathcal{O}(\epsilon)$ and collecting the terms that are independent of $\bs{p}$, yield
\begin{equation}
\label{eq:leading-order-invariant}
\bs{B}\bs{W}_{\mathbf{I}}\bs{S}_{\bs{0}}(\phi)+\Omega\bs{B}D_{\phi}\bs{X}_{\bs{0}}(\phi)=\bs{A}\bs{X}_{\bs{0}}(\phi)+\bs{F}^{\mathrm{ext}}(\phi),
\end{equation}
where $\bs{W}_\bs{I}=(\mathbf{W}_{ \mathbf{e}_1},\cdots,\mathbf{W}_{ \mathbf{e}_M})=(\bs{v}_1,\cdots,\bs{v}_M)$. Let $\bs{F}^\mathrm{ext}=\bs{F}^\mathrm{a}e^{\mathrm{i}\phi}+\bs{F}^\mathrm{a}e^{-\mathrm{i}\phi}$, and then we substitute the ansatz
\begin{gather}
\bs{X}_{\bs{0}}(\phi)=\bs{x}_{\bs{0}}e^{\mathrm{i}\phi}+\bar{\bs{x}}_{\bs{0}}e^{-\mathrm{i}\phi},\quad
\bs{S}_{\bs{0}}(\phi)=\bs{s}_{\bs{0}}^+e^{\mathrm{i}\phi}+{\bs{s}}_{\bs{0}}^-e^{-\mathrm{i}\phi},
\end{gather}
into~\eqref{eq:leading-order-invariant} and collect the coefficients of $e^{\mathrm{i}\phi}$ and $e^{-\mathrm{i}\phi}$, yielding
\begin{gather}
(\bs{A}-\mathrm{i}\Omega\bs{B})\bs{x}_{\bs{0}}=\bs{B}\bs{W}_{\mathbf{I}}\bs{s}_{\bs{0}}^+-\bs{F}^\mathrm{a}\label{eq:expphi},\\
(\bs{A}+\mathrm{i}\Omega\bs{B})\bar{\bs{x}}_{\bs{0}}=\bs{B}\bs{W}_{\mathbf{I}}{\bs{s}}_{\bs{0}}^--{\bs{F}}^{\mathrm{a}}\label{eq:expNegphi}.
\end{gather}
If $(\bs{A}-\mathrm{i}\Omega\bs{B})$ is nonsingular, we can simply set $\bs{s}_{\bs{0}}^+=\bs{0}$ and solve the linear system~\eqref{eq:expphi} to obtain $\bs{x}_{\bs{0}}$. However, if there exist eigenvalues near $\mathrm{i}\Omega$, e.g., $\lambda_i^{\mathcal{E}}\approx\mathrm{i}\Omega$, the coefficient matrix is nearly singular, we should choose $\bs{s}_{\bs{0}}$ such that the right-hand side vector is in the range of $(\bs{A}-\mathrm{i}\Omega\bs{B})$. This can be done by imposing an orthogonality constraint between the right-hand side vector and the kernel of $(\bs{A}-\mathrm{i}\Omega\bs{B})^\ast$~\cite{Jain2021HowModels}, i.e.,
\begin{equation}
(\bs{u}_i^\mathcal{E})^\ast\bs{B}\bs{W}_{\mathbf{I}}\bs{s}_{\bs{0}}^+-(\bs{u}_i^\mathcal{E})^\ast\bs{F}^\mathrm{a}=0.
\end{equation}
Once $\bs{s}_{\bs{0}}^+$ is obtained, we further solve for $\bs{x}_\bs{0}$. Likewise, we can solve for $\bs{s}_{\bs{0}}^-$.

\begin{remark}
\label{rk:ti-vs-tv}
We note that $\bs{x}_\bs{0}$ is $\Omega$ dependent and needs to be solved when $\Omega$ is changed. Therefore, solving for $\bs{x}_\bs{0}$ for a large collection of sampled $\Omega$ can be time-consuming for high-dimensional mechanical systems. In practice, one can check whether it is necessary to account for the contributions of $\bs{x}_\bs{0}$~\cite{li2024fast}. Specifically, for $\epsilon\ll1$, one may simply approximate the SSM expansion~\eqref{eq:ExpandW_e} as $\bs{z}(t)=\bs{W}(\bs{p}(t))$, which avoids the need to solve the linear systems~\eqref{eq:expphi} and~\eqref{eq:expNegphi}. Similar to the study~\cite{li2024fast}, we refer to the SSM solution as time-independent (TI) if we simply approximate the full system trajectory as $\bs{z}(t)\approx \bs{W}(\bs{p}(t))$. Furthermore, we refer to the SSM solution as time-varying (TV) if we include leading-order non-autonomous terms in its approximation, i.e., $\bs{z}(t)\approx \bs{W}(\bs{p}(t)) + \epsilon\bs{X}_{\bs{0}}(\Omega t) $. 
\end{remark}

\section{Non-intrusive SSM computation}
\label{sec:SemiIntrusive}

The main bottleneck of SSM computations is the composition $\bs{F}\circ\bs{W}$~\cite{Haro2016TheManifolds,Jain2021HowModels,thurnher2023nonautonomous}. Indeed, the evaluation of $\mathbf{F \circ W}\Bigr|_{\mathbf{m}}$ in the right-hand side of system~\eqref{app_eq:Hom_Aut_First} typically requires the expansion coefficients of the nonlinearity function $\bs{F}$, which renders SSM computations intrusive~\cite{Jain2021HowModels}.  These intrusive computations are inapplicable for generic FE packages because the function $\bs{F}$  is available only as a black-box input-output relation. Here, we develop expressions for polynomial composition $\mathbf{F \circ W}\Bigr|_{\mathbf{m}}$ based on the STEP~\cite{Muravyov2003DeterminationStructures}, thereby, facilitating non-intrusive SSM computations up to arbitrary polynomial orders.

\subsection{Split of odd and even nonlinearities}
We express the nonlinear vector-valued function $\bs{F}(\bs{z})$ in terms of its quadratic and cubic components as
\begin{align}\label{eq:NonintrusiveFunction}
    \mathbf{F}(\mathbf{z}) = \mathbf{F}_2 (\mathbf{z}) + \mathbf{F}_3 (\mathbf{z}).
\end{align}
We first demonstrate that the quadratic and cubic terms $\mathbf{F}_2$ and $\mathbf{F}_3$ can be computed non-intrusively simply via evaluations of the nonlinearity $\mathbf{F}$. Indeed, we note that $\mathbf{F}_2$ is an even function, and $\mathbf{F}_3$ is an odd function. Then, we can separate the two functions as
\begin{equation}
    \mathbf{F}_2(\mathbf{z}) = (\mathbf{F}(\mathbf{z}) + \mathbf{F}(-\mathbf{z}))/2,\quad
    \mathbf{F}_3(\mathbf{z}) = (\mathbf{F}(\mathbf{z}) - \mathbf{F}(-\mathbf{z}))/2.
\end{equation}
Based on the above separation, we have
\begin{equation}
\label{eq:comp}
    \mathbf{F \circ W}\Bigr|_{\mathbf{m}} = \mathbf{F}_2 \circ \bs{W}\Bigr|_{\mathbf{m}}+\mathbf{F}_3 \circ \bs{W}\Bigr|_{\mathbf{m}}.
\end{equation}
Next, we show that we can compute $\mathbf{F}_2 \circ \bs{W}$ and $\mathbf{F}_2 \circ \bs{W}$ in a non-intrusive fashion. For the purposes of illustration, we work with polynomial representations of the nonlinear functions as
\begin{align}
    \mathbf{F}_2 (\mathbf{z}) &=  \sum_{j,k=1}^N \mathbf{F}_{jk} {z}_j {z}_k,
     \ \ \ \quad \mathbf{F}_{jk} \in \mathbb{R}^{2n} , \\
    \mathbf{F}_3 (\mathbf{z}) &=  \sum_{j,k,l=1}^N\mathbf{F}_{jkl} {z}_j {z}_k{z}_l, \quad \mathbf{F}_{jkl} \in \mathbb{R}^{2n}.
\end{align}
Here and below, a vector's superscript denotes the vector's entry index. For instance, $\bs{z}^j$ represents the $j$-th entry of the vector $\bs{z}$.
As seen in~\eqref{eq:comp}, these functions need to be composed with the SSM parametrization $\mathbf{W}(\mathbf{p})$. In addition, the terms corresponding to distinct multi-index vectors need to be separated. Next, we derive the non-instrusive composition for the two nonlinear functions $\mathbf{F}_2$ and $\mathbf{F}_3$ via STEP as follows.

\subsection{Quadratic nonlinearity}
The composition of the quadratic nonlinearity $\mathbf{F}_2$ with the autonomous SSM parameterization $\mathbf{W}$ gives rise to pairs of multiplied parameterization coefficients
\begin{align}
    \mathbf{F}_2 \circ\mathbf{W} = \mathbf{F}_2 (\mathbf{W(p)}) &=  \sum_{j,k=1}^N \mathbf{F}_{jk} \left(\sum_{\mathbf{m}_1\in \mathbb{N}^M} \bs{W}_{\mathbf{m}_1,j} \mathbf{p}^{\mathbf{m}_1}\right)\left(\sum_{\mathbf{m}_2\in \mathbb{N}^M} \bs{W}_{\mathbf{m}_2,k} \mathbf{p}^{\mathbf{m}_2}\right)
    \nonumber\\
    &= \sum_{j,k=1}^N  \mathbf{F}_{jk} \sum_{\mathbf{m}_1,\mathbf{m}_2\in \mathbb{N}^M} \bs{W}_{\mathbf{m}_1,j} \bs{W}_{\mathbf{m}_2,k} \mathbf{p}^{\mathbf{m}_1+\mathbf{m}_2},
\end{align}
where $\bs{W}_{\mathbf{m}_1,j}$ stands for the $j$-th entry of the vector $\bs{W}_{\mathbf{m}_1}$ and $\bs{W}_{\mathbf{m}_2,k}$ denotes the $k$-th entry of the vector $\bs{W}_{\mathbf{m}_2}$.
Here the main challenge is to collect all terms for which $\mathbf{m}_1 + \mathbf{m}_2 = \mathbf{m}$ for a specific multi-index vector $\mathbf{m}$, and evaluate these terms simply via $\mathbf{F}_2$. In the special case that $\mathbf{m}_1=\mathbf{m}_2=\mathbf{m}/2$, we have
\begin{equation}
    \sum_{j,k=1}^N\mathbf{F}_{jk} \bs{W}_{\mathbf{m}_1,j} \bs{W}_{\mathbf{m}_1,k}
    = \mathbf{F}_2 (\mathbf{W}_{\mathbf{m}_1}).
\end{equation}
On the other hand, if $\mathbf{m}_1\neq\mathbf{m}_2$ but $\mathbf{m}_1 + \mathbf{m}_2 = \mathbf{m}$ still holds, the composed term can be obtained as
\begin{equation}
    \sum_{j,k=1}^N\mathbf{F}_{jk}  ( \mathbf{W}_{\mathbf{m}_1,j} \mathbf{W}_{\mathbf{m}_2,k} + \mathbf{W}_{\mathbf{m}_2,j} \mathbf{W}_{\mathbf{m}_1,k})
    = \frac{1}{2}\big(
    \mathbf{F}_2 (\mathbf{W}_{\mathbf{m}_1} + \mathbf{W}_{\mathbf{m}_2})
    - 
    \mathbf{F}_2 (\mathbf{W}_{\mathbf{m}_1} - \mathbf{W}_{\mathbf{m}_2})\big).
\end{equation}
Therefore, we can indeed evaluate $[\ \mathbf{F}_2 \circ \mathbf{W} ]_{\mathbf{m}}$ in a non-intrusive manner.

\subsection{Cubic nonlinearity}
The composition of the cubic nonlinearity $\mathbf{F}_3$ with the autonomous SSM parameterization $\mathbf{W}$ gives rise to pairs of multiplied parameterization coefficients as
\begin{align}
    \bs{F}_3\circ\bs{W}=\mathbf{F}_3 (\mathbf{W(p)}) &=  \sum_{j,k,l=1}^N \mathbf{F}_{jkl} \left(\sum_{\mathbf{m}_1\in \mathbb{N}^M} \bs{W}_{\mathbf{m}_1,j} \mathbf{p}^{\mathbf{m}_1}\right)\left(\sum_{\mathbf{m}_2\in \mathbb{N}^M} \bs{W}_{\mathbf{m}_2,k} \mathbf{p}^{\mathbf{m}_2}\right) \left(\sum_{\mathbf{m}_3\in \mathbb{N}^M} \bs{W}_{\mathbf{m}_3,l} \mathbf{p}^{\mathbf{m}_3}\right)
    \nonumber\\
    &= \sum_{j,k,l=1}^N  \mathbf{F}_{jkl} \sum_{\mathbf{m}_1,\mathbf{m}_2,\mathbf{m}_3\in \mathbb{N}^M} \bs{W}_{\mathbf{m}_1,j} \bs{W}_{\mathbf{m}_2,k} \bs{W}_{\mathbf{m}_3,l} \mathbf{p}^{\mathbf{m}_1+\mathbf{m}_2+\mathbf{m}_3}.
\end{align}
Here $\bs{W}_{\mathbf{m}_3,l}$ stands for the $l$-th entry of the vector $\bs{W}_{\mathbf{m}_3}$. Similarly, the main challenge is to collect all terms for which $\mathbf{m}_1 + \mathbf{m}_2 + \mathbf{m}_3 = \mathbf{m}$ for a specific multi-index vector $\mathbf{m}$, and evaluate these terms simply via $\mathbf{F}_3$.  For brevity we introduce the following notation $\mathbf{v}_i := \mathbf{W}_{\mathbf{m}_i}$. In the special case that $\mathbf{m}_1=\mathbf{m}_2=\mathbf{m}_3=\mathbf{m}/3$, we have
\begin{equation}
   \sum_{j,k,l=1}^N \mathbf{F}_{jkl} \mathbf{v}_{1,j} \mathbf{v}_{1,k} \mathbf{v}_{1,l} = \mathbf{F}_3 (\mathbf{v}_1).
\end{equation}
In the special case that two of the multi-index vectors are the same, we obtain
\begin{equation}
  \sum_{j,k,l=1}^N \mathbf{F}_{jkl}  \ \ \ \mathclap{ \sum_{(a,b,c) \in \sigma[1,1,2]}} \quad \mathbf{v}_{a,j} \mathbf{v}_{b,k} \mathbf{v}_{c,l}
    =
    \bigg(\mathbf{F}_3 (\mathbf{v}_1+\mathbf{v}_{2}) 
    - \mathbf{F}_3 (\mathbf{v}_{1} -\mathbf{v}_{2})  
    - 2 \mathbf{F}_3 (\mathbf{v}_{2}) \bigg)/2,
\end{equation}
where $\sigma[1,1,2]=\{(1,1,2),(1,2,1),(2,1,1)\}$, denoting the collection of permutations of the vector $[1,1,2]$. Finally, if all three multi-index vectors are distinct, we have
\begin{align} \nonumber
   \sum_{j,k,l=1}^N \mathbf{F}_{jkl}   \ \ \ \ \mathclap{ \sum_{(a,b,c) \in \sigma[1,2,3]}} \quad \mathbf{v}_{a,j} \mathbf{v}_{b,k} \mathbf{v}_{c,l}
    =
    &\mathbf{F}_3 (\mathbf{v}_{1}+\mathbf{v}_{2}+ \mathbf{v}_{3}) 
    -\mathbf{F}_3 (\mathbf{v}_{1}+ \mathbf{v}_{2}) 
    \\
    &-\mathbf{F}_3 (\mathbf{v}_{1}+\mathbf{v}_{3}) 
    -\mathbf{F}_3 (\mathbf{v}_{2}+ \mathbf{v}_{3}) 
    +\sum_{i=1}^3 \mathbf{F}_3 (\mathbf{v}_{{i}}) .
\end{align}
where $\sigma[1,2,3]$ denotes the collection of permutations of the vector $[1,2,3]$. Therefore, we can evaluate $[\bs{F}_3\circ\bs{W}]_\bs{m}$ in a non-intrusive manner as well. This allows us to evaluate $[\bs{F}\circ\bs{W}]_\bs{m}$ and solve system~\eqref{app_eq:Hom_Aut_First} for each $\mathbf{m}$ non-intrusively.

\section{Implementation}
\label{sec:implementation}
We have implemented the above non-intrusive algorithm in SSMTool~\cite{SSMTool2}, an open-source package for the automated computation of SSMs and their associated ROMs. To demonstrate the effectiveness of the developed non-intrusive algorithm, we consider two FE software: an open-source MATLAB-based package, YetAnotherFECode~\cite{YetAnotherFEcode}, and a commercial package, COMSOL Multiphysics. We couple SSMTool with these packages while treating them as black boxes for the nonlinear function handle $\bs{F}$ to perform SSM-based model reduction of FE models. In particular, FE models are constructed using the FE toolboxes. Then communications between SSMTool and FE toolboxes will be built such that the FE toolboxes will return the output $\bs{F}(\bs{u})$ when an input vector $\bs{u}$ is specified by SSMTool. We note that the communication time between SSMTool and the in-house YetAnotherFECode is minimal because the package has been integrated well into SSMTool. On the other hand, such communication time is significant for the commercial software COMSOL Multiphysics. We also note that the coupling established here can be easily generalized to other commercial or open-source FE code.

%\subsection*{Treatment for complex inputs}
As seen in the previous section, the two nonlinear function handles $\mathbf{F}_2$ and $\mathbf{F}_3$ will take the SSM parameterization $\mathbf{W}_{\mathbf{m}_i}$ as inputs. We note that these inputs are typically complex-valued vectors. However, we find that the function handles optimized for function evaluations over real-valued inputs. Indeed, incorrect evaluations will be returned if complex vectors are taken as inputs in commercial packages such as COMSOL 
Multiphysics. To resolve this issue, we use a combination of function evaluations on a complex input's real and imaginary parts to reconstruct the complex nonlinearity. In particular, for quadratic nonlinearities, consider a complex vector $\mathbf{v}$, define
\begin{equation}
\mathbf{v}_a=\mathrm{Re}(\mathbf{v}),\quad \mathbf{v}_b=\mathrm{Im}(\mathbf{v}),\quad \mathbf{v}'=\mathbf{v}_a+\mathbf{v}_b,
\end{equation}
we obtain
\begin{equation}
\mathbf{F}_2(\mathbf{v}) = \mathrm{i} \mathbf{F}_2(\mathbf{v}')+(1-\mathrm{i})\mathbf{F}_2(\mathbf{v}_a)-(1+\mathrm{i})\mathbf{F}_2(\mathbf{v}_b).
\end{equation}
Likewise, for cubic nonlinearities, we further define $\mathbf{v}''=\mathbf{v}_a-\mathbf{v}_b$, and we have
\begin{equation}
\mathbf{F}_3(\mathbf{v}) = 2\mathbf{F}_3(\mathbf{v}_a)+\frac{-1+\mathrm{i}}{2}\mathbf{F}_3(\mathbf{v}')-\frac{1+\mathrm{i}}{2}\mathbf{F}_3(\mathbf{v}'')-2\mathrm{i}\mathbf{F}_3(\mathbf{v}_b).
\end{equation}

\section{Examples}
\label{sec:examples}
We consider a suite of examples to demonstrate the effectiveness of the non-intrusive algorithm for the computation of SSMs and their associated ROMs. The mass and stiffness matrices are positive definite in all FE examples reported below, and $\bs{C}$ is a Rayleigh damping matrix. However, as stated, our non-intrusive implementation works for general mechanical systems. To illustrate its wide applicability, we add a cantilevered pipe conveying fluid example in~\ref{sec:add-example}. This pipe system has asymmetric damping and stiffness matrices resulting from flow-induced gyroscopic and follower forces and velocity-dependent nonlinear internal forces~\cite{li2023model}. As shown in~\ref{sec:add-example}, our non-intrusive procedure indeed works for such a general mechanical system.

\subsection{A shallow arc with 1:2 internal resonances}

\begin{figure}[!ht]
\centering
\includegraphics[width=0.4\textwidth]{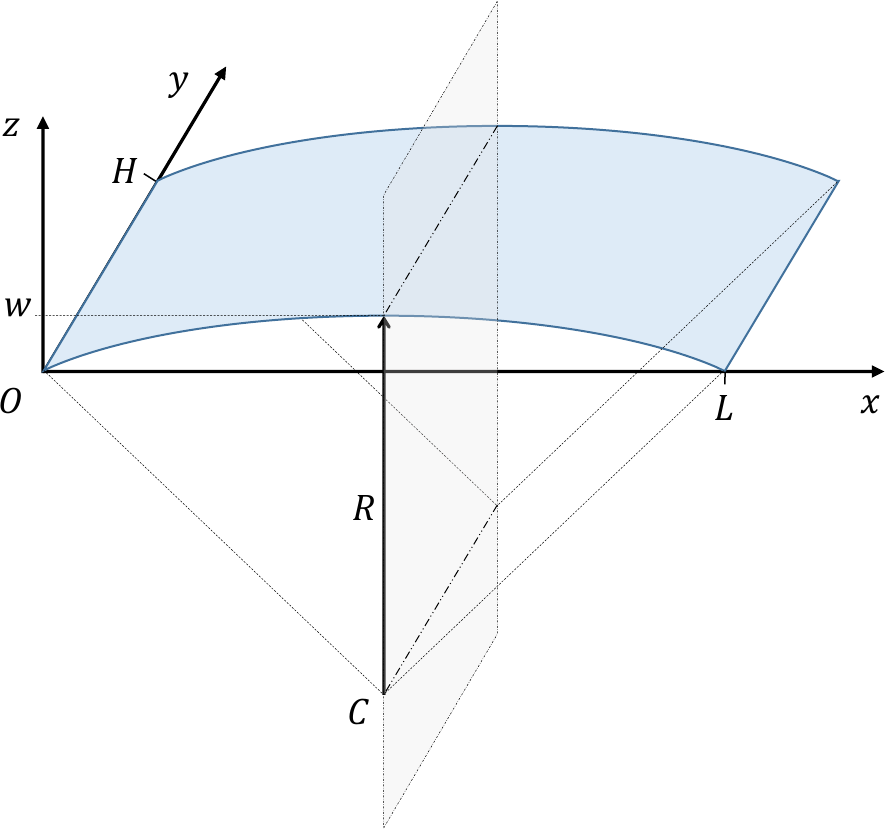}
\caption{The schematic of a shallow shell structure~\cite{Jain2021HowModels,part-i}.}
\label{fig:shell_mesh}
\end{figure}

As our first example, we consider the nonlinear vibrations of a shallow-arc structure~\cite{part-i}, shown in Fig.~\ref{fig:shell_mesh}. Here, the shell is simply supported at the two opposite edges aligned along the $y-$axis in Fig.~\ref{fig:shell_mesh}. We refer to~\cite{part-i} for this shell's geometric and material properties. We adopt the same finite-element model as in~\cite{part-i}, and the discrete model has 400 elements and 1,320 DOFs. As demonstrated in~\cite{part-i}, one can tune the curvature of this shell such that the first two bending modes of this structure admit a 1:2 internal resonance~\cite{part-i}. Specifically, the eigenvalues of the first two pairs of modes of the FE model are listed below~\cite{part-i}
\begin{equation}
    \lambda_{1,2}=-0.30\pm\mathrm{i}149.22,\quad\lambda_{3,4}=-0.60\pm\mathrm{i}298.78.
\end{equation}

\begin{figure}[!ht]
\centering
\includegraphics[width=0.6\textwidth]{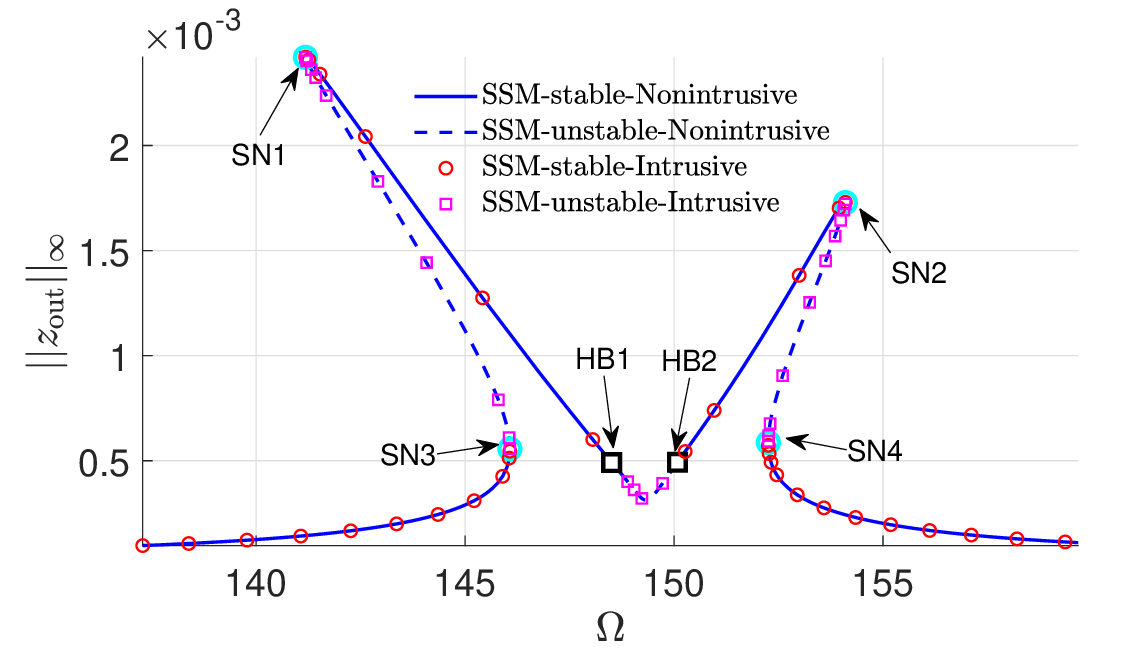}
\caption{The forced response curve for the shallow shell structure with 1:2 internal resonance. Here, the solid and dashed lines denote stable and unstable periodic orbits predicted using non-intrusive SSM reduction, the red circles and magenta squares represent stable and unstable periodic orbits obtained from intrusive SSM reduction, the cyan circles denote saddle-node (SN) bifurcation points, and black squares denote Hopf bifurcation (HB) points. All results here are obtained from $\mathcal{O}(5)$ approximations of the SSM.}
\label{fig:shell_FRC}
\end{figure}

A concentrated load $\epsilon100\cos\Omega t$ in the $z-$ direction is applied at mesh node with $(x,y)=(0.25L,0.5H)$. We are concerned with the forced response around the first mode regarding the $z$-displacement of the node. We perform reduction on four-dimensional SSM to account for the internal resonance. Thus, the dimension of the phase space is reduced from 2640 to four. As demonstrated in~\cite{li2024fast}, the TV-SSM solution is in good agreement with the TI-SSM solution, and hence, we adopt the TI-SSM solution in this example (see Remark~\ref{rk:ti-vs-tv} for more details on these two types of solutions). As illustrated in~\cite{part-i}, the forced response curve of this system converges well at $\mathcal{O}(5)$ approximation. We plot the FRC obtained using the $\mathcal{O}(5)$ approximation at Fig.~\ref{fig:shell_FRC}. In particular, we present the FRC obtained from the intrusive computation~\cite{part-i} and the non-intrusive scheme proposed here. As shown in Fig.~\ref{fig:shell_FRC}, the FRCs obtained from the two computations match perfectly, which validates the effectiveness of the non-intrusive computations.

From the FRC in Fig.~\ref{fig:shell_FRC}, we observe four bifurcated periodic orbits along the curve. Two are saddle-node bifurcations, and the other are secondary Hopf (torus) bifurcations.  We again see the convergence of the critical excitation frequencies $\Omega_{\mathrm{SN1}}$ and $\Omega_{\mathrm{HB1}}$ as the expansion order is increased. 

% Here, we also list the critical excitation frequency $\Omega$, at which the first saddle-node bifurcation (SN1) and Hopf bifurcation (HB1) occur, in Table.~\ref{omega-bif-shell}. Moreover, an excellent match regarding these critical excitation frequencies between the intrusive and non-intrusive computations is again obtained because their relative errors are consistently less than $10^{-9}$. This again demonstrates the accuracy of the non-intrusive computations.

% \begin{table*}[ht]
% \begin{center}
% \begin{tabular}{ c|c|c }
% \hline
% & $\mathcal{O}(3)$ & $\mathcal{O}(5)$ \\
% \hline
% $\Omega_\mathrm{SN1,Intrusive}$ & 141.5837356784178 & 141.1805294942144 \\
% \hline
% $\Omega_\mathrm{SN1,Non-intrusive}$ & 141.5837356780345 & 141.1805294963087 \\
% \hline
% $\Omega_\mathrm{HB1,Intrusive}$ & 148.5069797491482 & 148.5064945762714 \\
% \hline
% $\Omega_\mathrm{HB1,Non-intrusive}$ & 148.5069803419740 & 148.5064946216117 \\
% \hline
% \end{tabular}
% \end{center}
% \caption{Critical excitation frequencies for the appearance of bifurcation points SN1 and HB1 shown in Fig.~\ref{fig:shell_FRC}. These bifurcation points are persistent and converged as the expansion order is increased from $\mathcal{O}(3)$ to $\mathcal{O}(5)$. Here, $\Omega_\mathrm{bif,comp}$ with $\mathrm{bif}\in\{\mathrm{SN1},\mathrm{HB1}\}$ denoting the bifurcation points and $\mathrm{comp}\in\{\mathrm{Intrusive},\mathrm{Non-intrusive}\}$ providing the computational schemes.}
% \label{omega-bif-shell}
% \end{table*}

In this example, the estimated memory cost for the non-intrusive and intrusive algorithms at $\mathcal{O}(5)$ computation are around 8 Mb and 35 Mb. The latter is more than 5 times that of the former. This demonstrates that the non-intrusive algorithm can indeed reduce memory consumption. The computational time of the non-intrusive and intrusive algorithms are about 6 minutes and one minute. This suggests that the implementation of the nonlinear function handle is not efficient. Indeed, during the 6 minutes, nearly 97\%  of them were used by the evaluations of the nonlinear function handles. Thus, one can significantly speed up the computations by improving the computational efficiency of the nonlinear function handles.

\subsection{A NACA wing}

As our second example, we consider a geometrically nonlinear FE model of an aircraft-wing~\cite{Jain2017ADynamics} with 49968 triagonal shell elements and, accordingly, 133920 DOFs. Model reduction has been performed for this high-dimensional mechanics problem via quadratic manifolds~\cite{Jain2017ADynamics} and also spectral submanifolds~\cite{Jain2017ADynamics}. In particular, intrusive computations were used in~\cite{Jain2017ADynamics} to construct a two-dimensional SSM, and the associated ROM enables analytic prediction on the forced response curve. Here, we apply our proposed non-intrusive algorithm to perform the same reduction and compare the results obtained from the intrusive and non-intrusive computations to validate the effectiveness of the non-intrusive algorithm.

\begin{figure}[!h]
\centering
\includegraphics[width=.4\textwidth]
{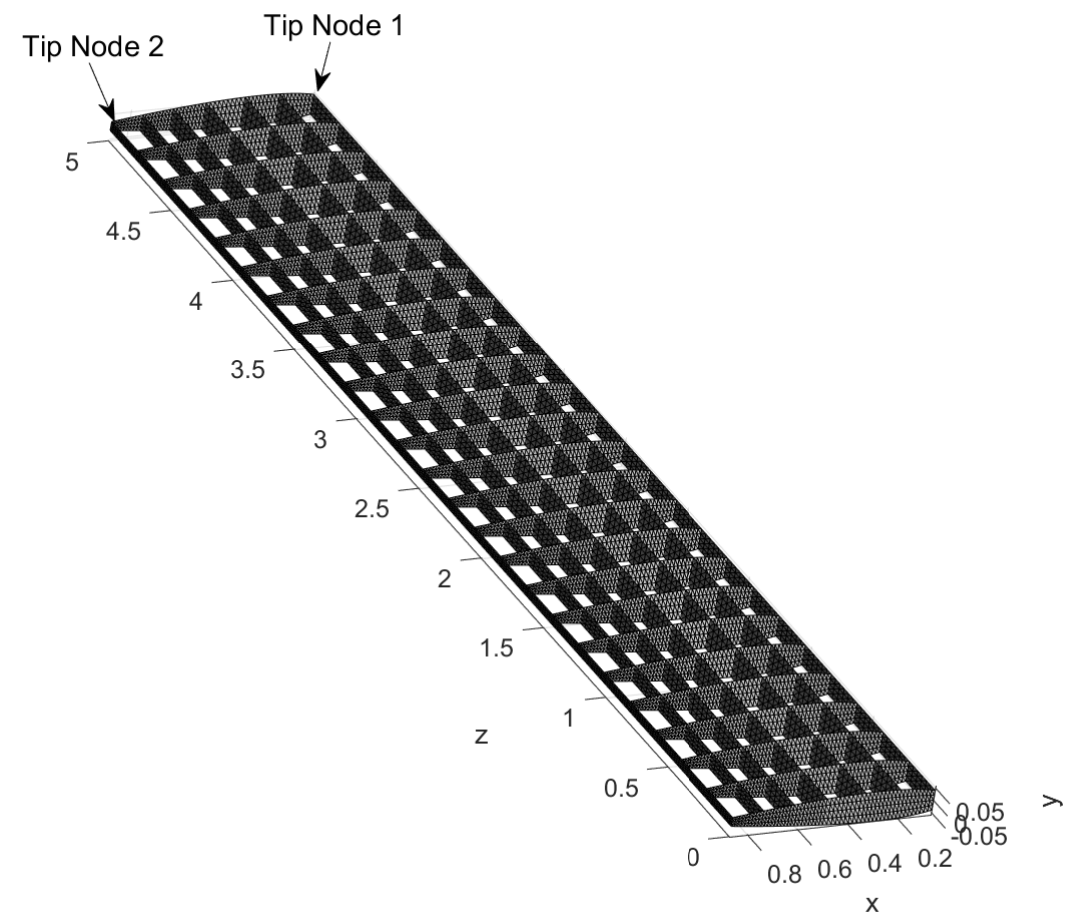}
\includegraphics[width=.45\textwidth]
{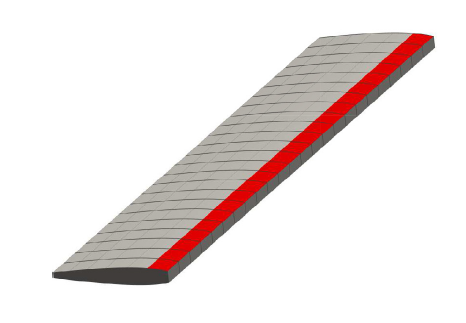}
\caption{(Left panel) Model showing the interior structure of the wing. These ribs are then covered with plates to create the model considered here \cite{Jain2021HowModels}. (Right panel) The full model of the wing is covered with plates. Boundary conditions are chosen, such as fixing the end of the wing as if attached to an aircraft.}
\label{fig:Wing}
\end{figure}

A depiction of the supporting structure of the wing, as well as the model of the wing, is shown in Fig.~\ref{fig:Wing}. Clamped boundary conditions at $z=0$ are applied as if it is attached to an aircraft. The values of physical parameters of this wing are chosen according to Table 7 of~\cite{Jain2021HowModels}. We adopt Rayleigh damping with damping constants chosen such that the overall damping ratio is $0.4\%$~\cite{Jain2021HowModels}, which produces a strongly underdamped structure. In particular, the first pair of eigenvalues is given as $\lambda_{1/2} = -0.058 \pm  29.343 \mathrm{i}$.

\begin{figure}[!ht]
\centering
\includegraphics[width=.45\textwidth]{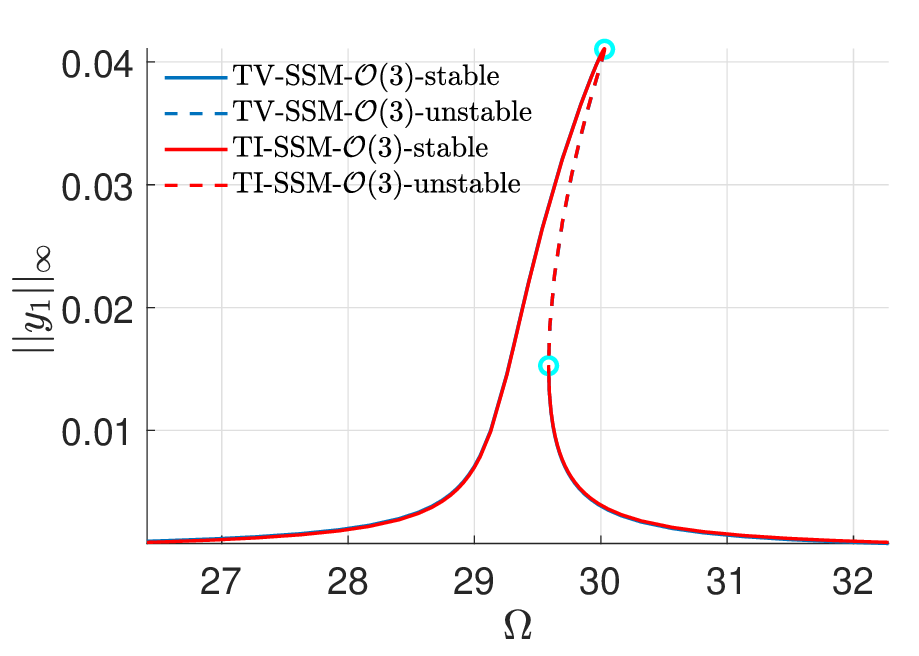}
\includegraphics[width=.45\textwidth]{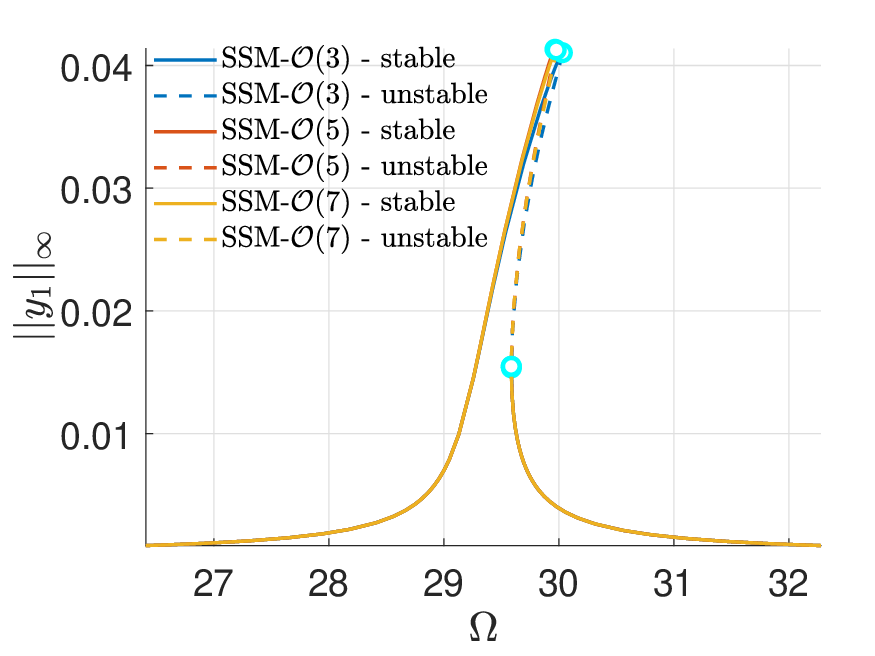}\\
\includegraphics[width=.45\textwidth]{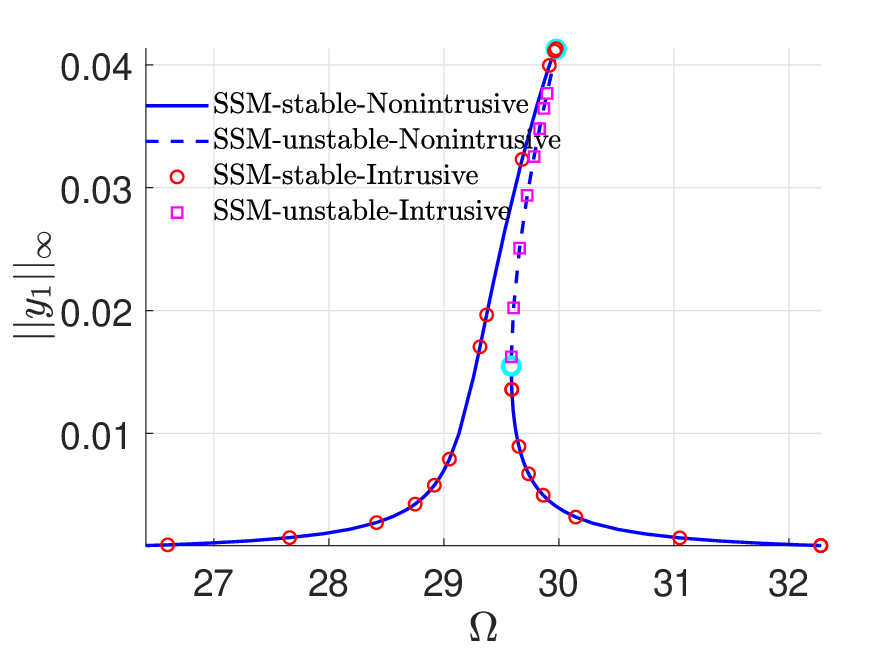}
\caption{Forced response curve regarding the tip deflection of the aircraft wing. In the upper-left panel, we present the results obtained from $\mathcal{O}(3)$ approximation with and without the contribution of the leading-order non-autonomous part of SSM. Specifically, TV and TI stand for time-varying and time-independent solutions with and without the contribution. The results predicted from SSM-based reductions truncated at various expansion orders are shown in the upper-right panel. The lower panel presents the results at $\mathcal{O}(7)$ expansion, obtained using both the intrusive computation in~\cite{li2023nonlinear} and our proposed non-intrusive algorithm.}
\label{fig:wing_frcs}
\end{figure}

Following~\cite{Jain2021HowModels}, we apply two concentric harmonic excitation loads $10\cos\Omega t$ at the two tip nodes shown in Fig.~\ref{fig:Wing}. We aim to extract the FRC for the primary resonance of the first mode. We compute the FRC of the system based on TV-SSM and TI-SSM solutions to check whether TI-SSM is accurate enough. As seen in the upper-left panel of Fig.~\ref{fig:wing_frcs}, the FRC based on the TI-SSM solution matches that of the TV-SSM solution. Therefore, we conclude that the TI-SSM based predictions have sufficient accuracy and we will use the TI-SSM solution in the rest of this example.

Next, we compute the FRC with increasing expansion orders using both the intrusive algorithm~\cite{Jain2021HowModels} and the non-intrusive algorithm established here. As seen in the upper-right panel of Fig.~\ref{fig:wing_frcs}, the FRC converges well at $\mathcal{O}(5)$ approximation. Finally, the lower panel of Fig.~\ref{fig:wing_frcs} shows that the results from the non-intrusive algorithm are the same as that of the intrusive computation, as expected.

In this example, the estimated memory cost for the non-intrusive and intrusive algorithms at $\mathcal{O}(7)$ computation are around 335 Mb and 1182 Mb. The latter is nearly 3.5 times that of the former. This demonstrates that the non-intrusive algorithm can indeed reduce memory consumption. The computational time of the non-intrusive and intrusive algorithms is about 4 hours and 0.5 hours, respectively. This again suggests that implementing the nonlinear function handle is inefficient. Indeed, during the 4 hours, nearly 93\% of them were used for the evaluations of the nonlinear function handles. Thus, one can significantly speed up the computations by improving the computational efficiency of the nonlinear function handles.

\subsection{A viscoelastic perforated cover plate}

As our third example, we consider the FE model of a perforated cover plate~\cite{ehrhardt2017finite}, which is a part of the exhaust system of a large diesel engine. The welded boundary of the plate is approximated by a series of 80 linear springs in the radial direction, each having a stiffness of 650 $\mathrm{N/m}$. The diameter and thickness of the plate are 317.5 mm and 1.5 mm. Other geometric parameters, such as curvature, can be found in~\cite{ehrhardt2017finite}. Material parameters of the plate are listed below: a Young's modulus of 96 GPa, a Poisson ratio of 0.3, and a density of 5120 $\mathrm{kg/m^3}$. Here, we use a Kelvin-Voigt model to account for viscoelasticity. In particular, the viscosity parameter is chosen as 0.6 Mpa. Thus, this model has both displacement-dependent and velocity-dependent nonlinear internal forces.

\begin{figure}[!ht]
\centering
\includegraphics[width=0.45\textwidth]{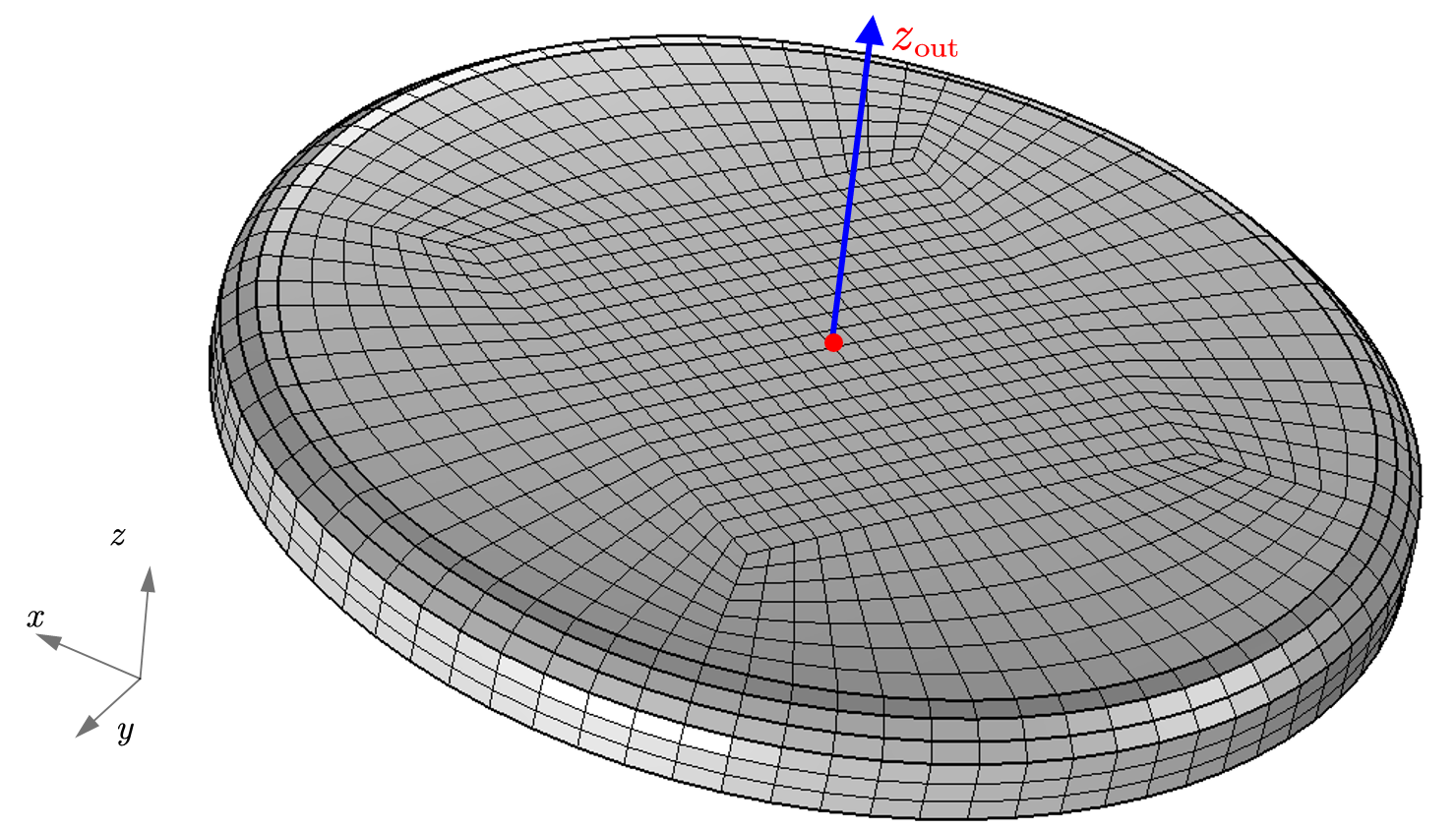}
\caption{A schematic plot of the mesh for FE model of the perforated cover plate.}
\label{fig:rom_plate_plot}
\end{figure}

The plate is discretized with 1440 quadrilateral elements, 400 edge elements, and 21 vertex elements, resulting in a FE model with 9846 DOFs. With boundary conditions imposed, the number of DOFs is decreased to 7405. Thus, the phase space of the full system here is of dimension 14810. The first two natural frequencies of the system are $\omega_1=1272.0~\mathrm{rad/s}$ and $\omega_2=2065.8~\mathrm{rad/s}$. Here, we take the two-dimensional subspace corresponding to the first mode as the master subspace for model reduction. 

\begin{figure}[!ht]
\centering
\includegraphics[width=0.45\textwidth]{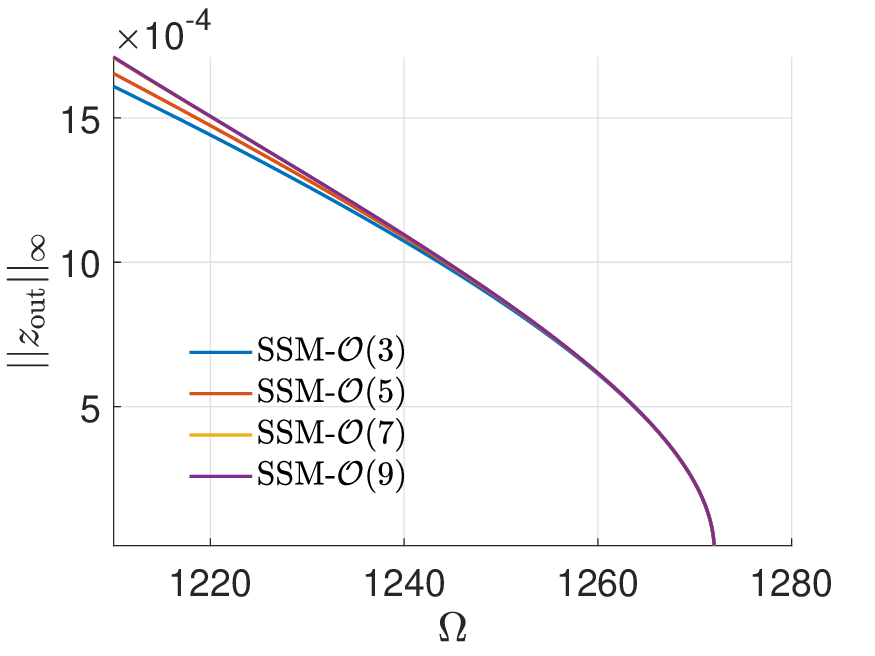}
\caption{Backbone curve for the first vibration mode of the perforated plate. The curve is obtained using SSM-based reduction truncated at various orders.}
\label{fig:rom-bb}
\end{figure}

We use the proposed non-intrusive algorithm to compute the autonomous SSM of the master subspace and then extract the backbone curve regarding the $z$-displacement at the center of the plate (cf. Fig.~\ref{fig:rom_plate_plot}). The obtained backbone curve is shown in Fig.~\ref{fig:rom-bb}, displaying a softening behavior. We observe that the $\mathcal{O}(3)$ truncation gives an accurate prediction when the displacement is less than half the thickness of the plate. In addition, an $\mathcal{O}(7)$ or higher-order approximation is needed to yield converged solutions when the displacement is larger than the thickness of the plate.

\begin{figure}[!ht]
\centering
\includegraphics[width=.45\textwidth]{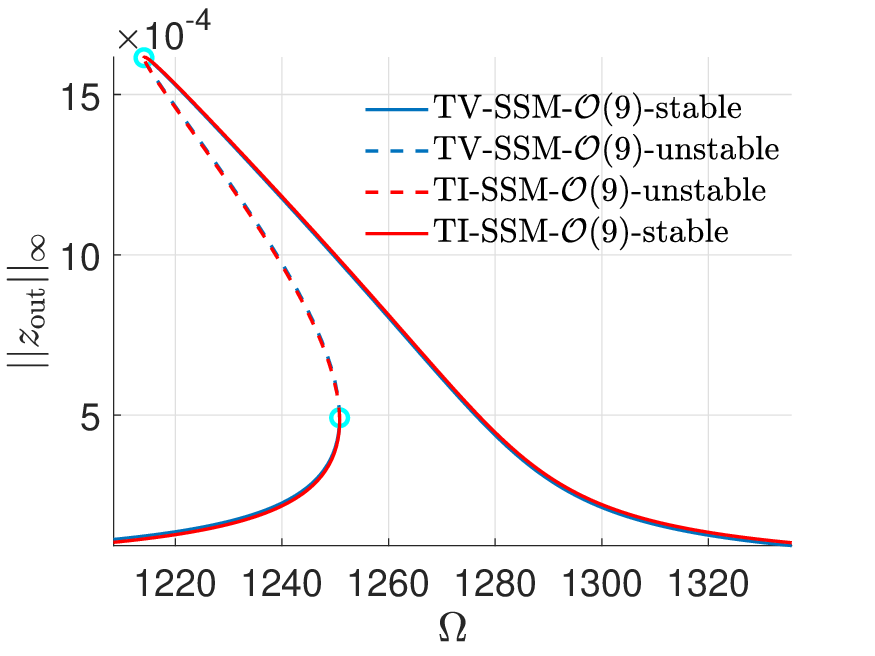}
\includegraphics[width=.45\textwidth]{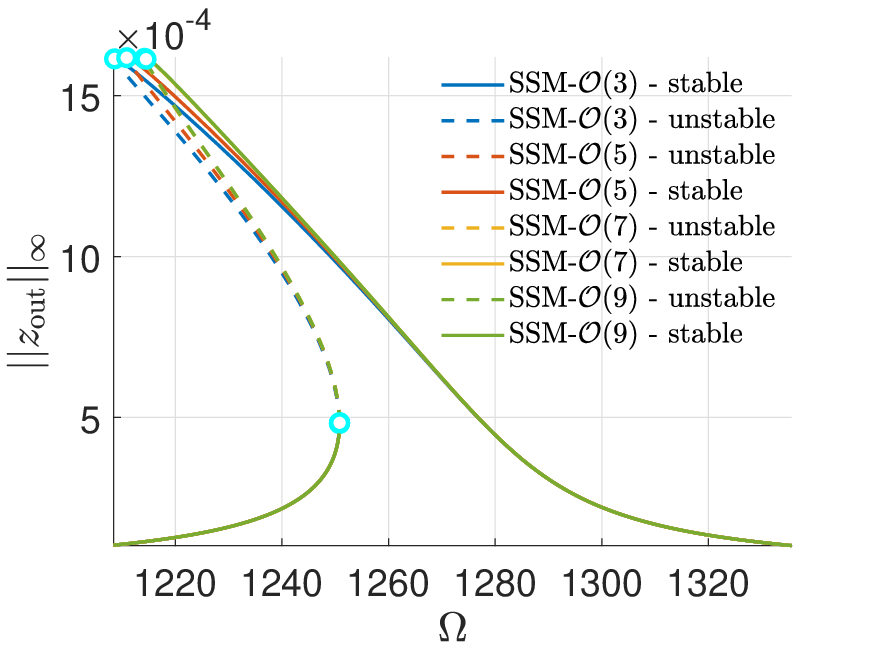}\\
\includegraphics[width=.45\textwidth]{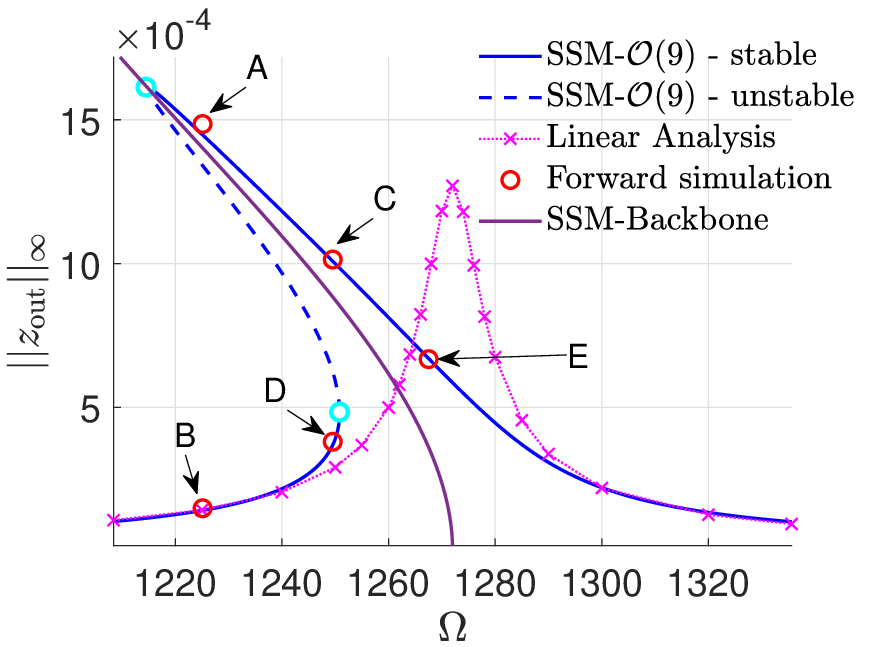}
\includegraphics[width=.45\textwidth]{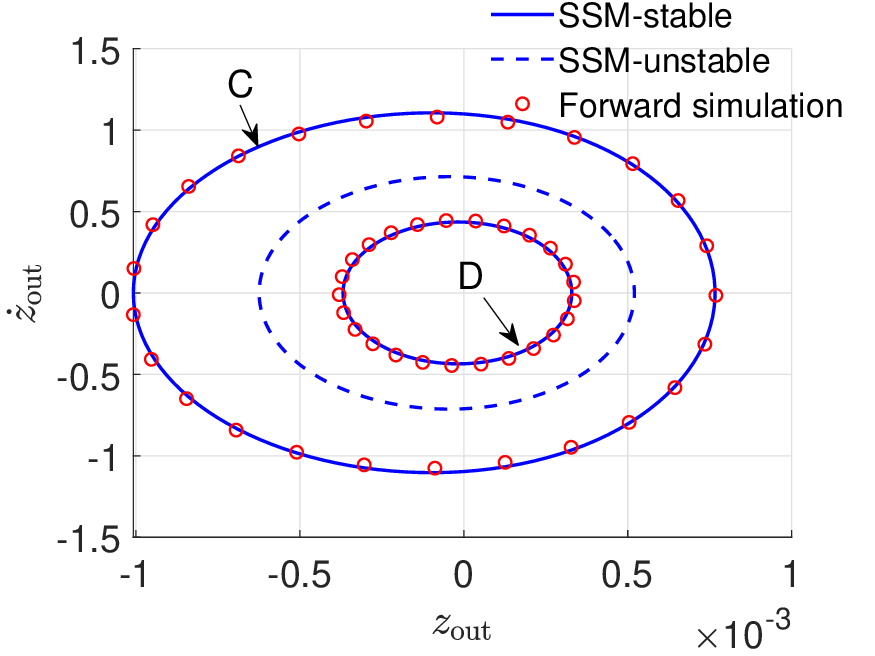}
\caption{Forced response regarding the deflection at the center of the perforated cover plate. In the upper-left panel, we present the forced response curve (FRC) obtained from $\mathcal{O}(3)$ approximation with and without the contribution of the leading-order non-autonomous part of SSM. Specifically, TV (TI) stands for time-varying (time-independent) solution with (without) the contribution. In the upper-right panel, the FRCs predicted from the SSM-based reduction truncated at various expansion orders are shown. The lower-left panel presents the FRC and the backbone curve at $\mathcal{O}(9)$ expansion, along with results of liner analysis and forward simulation applied to the full system. In the lower-right panel, the projection of the coexisting three periodic orbits at $\Omega=\Omega_\mathrm{C}=1249.5~\mathrm{rad/s}$ onto the plane $(z_\mathrm{out},\dot{z}_\mathrm{out})$ are shown (cf.~the lower-left panel). Here, the solid and dashed lines denote stable and unstable periodic orbits.}
\label{fig:rom_frcs}
\end{figure}

Next, we add a harmonic excitation $1.6\cos\Omega t$ at the center point along the $z$-axis (see Fig.~\ref{fig:rom_plate_plot}) and consider the primary resonance of the first mode of the perforated cover plate, namely, $\Omega\approx\omega_1$. Similar to the previous examples, we first check whether we can simply take TI-SSM solutions. As seen in the upper-left panel of Fig.~\ref{fig:rom_frcs}, the FRC using TV-SSM solutions matches well with that of the TI-SSM solutions. Therefore, it is sufficient to use TI-SSM solutions here. We then check the convergence of the FRC with increasing expansion orders. We observe from the upper-right panel of Fig.~\ref{fig:rom_frcs} the FRC is well converged at $\mathcal{O}(9)$ expansion because the FRC at $\mathcal{O}(7)$ coincides with the FRC at $\mathcal{O}(9)$.

We plot the converged FRC along with the backbone curve in the lower-left panel of Fig.~\ref{fig:rom_frcs}, from which we see the backbone curve passes through a saddle-node (SN) bifurcation point on the FRC. This SN point also marks the peak of the FRC. To validate the effectiveness of the predicted FRC, we first extract the linear response function of the full system. As seen in the figure, the SSM-based prediction matches well with the linear analysis when $\Omega$ is away from the resonance. However, significant discrepancies are observed near the resonance where the periodic response has a large amplitude, and geometric nonlinearity plays an important role. We carry out forward simulations of the full system to further validate the SSM-based nonlinear predictions. Specifically, we take five sampled points on the FRC, namely, points A-E in the lower-left panel of Fig.~\ref{fig:rom_frcs}. Here $\Omega_\mathrm{A}=\Omega_\mathrm{B}<\Omega_\mathrm{C}=\Omega_\mathrm{D}<\Omega_\mathrm{E}$. In each case, we initialize the forward simulation with a starting point on the periodic orbit predicted from the SSM-based ROM and perform the simulation with 100 cycles such that a steady state has been reached. We then extract the magnitude of these periodic responses at steady state and plot them in the lower-left panel of Fig.~\ref{fig:rom_frcs}, from which we see that the results from the forward simulation match well with the SSM-based predictions.

We note that there are three coexisting periodic orbits when $\Omega=\Omega_\mathrm{A}$ or $\Omega=\Omega_\mathrm{C}$. In fact, this coexistence holds if $\Omega\in(\Omega_\mathrm{SN1},\Omega_\mathrm{SN2})$, where $\Omega_\mathrm{SN1}$ and $\Omega_\mathrm{SN2}$ mark the critical forcing frequencies of the two SN bifurcation points. To have a close look at the coexistence, we plot the three periodic orbits at $\Omega=\Omega_\mathrm{C}=1249.5~\mathrm{rad/s}$ in the phase plane $(z_\mathrm{out},\dot{z}_\mathrm{out})$ at the lower-right panel of Fig.~\ref{fig:rom_frcs}, where we again see the SSM-based predictions match well with that of the forward simulations.

In this example, the estimated memory cost for the non-intrusive algorithm at $\mathcal{O}(9)$ computation is 26 Mb, and the associated computational time is 4.5 hours. In contrast, the averaged simulation time for the five forward simulations is more than 43 hours. This illustrates that a significant speed-up gain can be obtained using the SSM-based reduction. Indeed, we can even analytically extract the forced response surface under variations in forcing frequency and amplitude once the autonomous SSM is computed~\cite{li2024fast}. The forward simulations, in contrast, need to be repeated when the forcing parameters are changed and can only capture stable periodic orbits even with good initial guesses.

\subsection{A MEMS resonator}

\begin{figure}[!ht]
\centering
\includegraphics[width=0.45\textwidth]{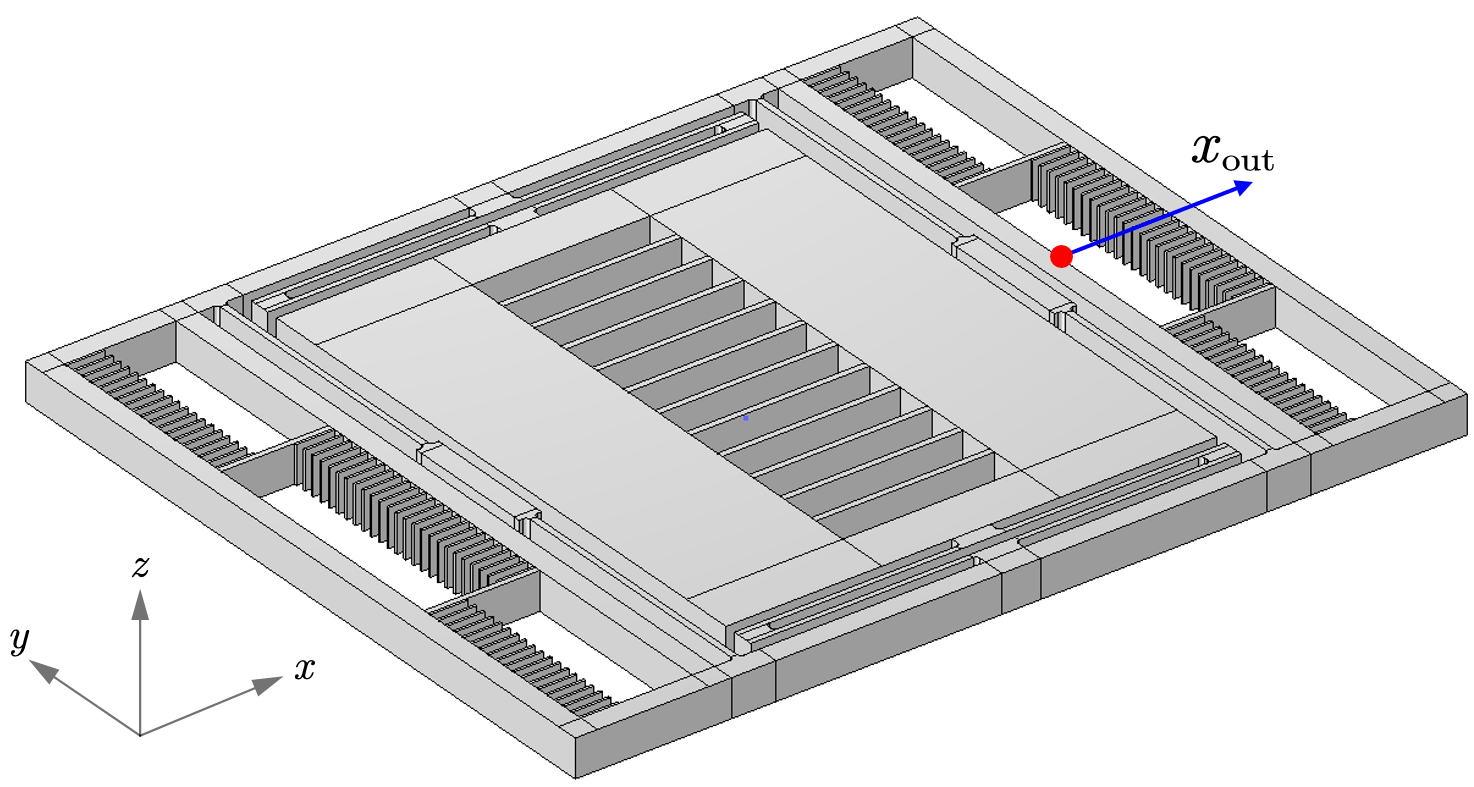}
\caption{A schematic plot of a MEMS gyroscope prototype.}
\label{fig:memsplot}
\end{figure}

As our last example, we consider a large FE model representing a MEMS gyroscope prototype~\cite{marconi2021exploiting}. As seen in Fig.~\ref{fig:memsplot}, this MEMS device consists of a frame and a proof mass. The frame is attached to the ground using four flexible beams, and the proof mass is connected to the frame using two flexible folded beams. In the first vibration mode, the frame and the proof mass oscillate synchronously along the $x$ direction. As for the second vibration mode, the frame stays nearly still, and the proof mass oscillates in the direction of $y$. Therefore, when the device is subject to external angular rotation along the $z$ axis, oscillating along the first mode will generate a Coriolis force on the proof mass along the $ y$ direction, which excites the second vibration mode. The vibration along the second mode is detected and further converted into an angular velocity measurement. This prototype exhibits a strongly nonlinear forced response along the first mode by design. Then, the drive frequency near the first mode can be tuned to match the sense frequency along the second mode~\cite{marconi2021exploiting}.

The overall physical dimensions of the device are $600\times600\times20~[\mu\mathrm{m}]$ with more than a million degrees of freedom. Detailed geometric and physical parameters of the structure can be found in~\cite{marconi2021exploiting}. The FE model consists of 28,084 hexahedral elements and 8,636 prismatic (wedge) elements, resulting in a high-dimensional discrete model with 1,029,456 DOFs. Here, we adopt a Rayleigh damping $\bs{C}=1659.8\bs{M}+1.336\times10^{-7}\bs{K}$ such that the first mode has a quality factor $Q=100$. In practice, the quality factor is much higher. Here, we consider this moderately low damping for the convenience of validation. In particular, we will use forward simulation of the full system to validate the predicted periodic orbits. We take the moderately low damping so that the forward simulations can reach steady states in a reasonable time. We note that our non-intrusive algorithm, however, also works well in the case of large quality factors. Moreover, the computational cost of constructing the SSM-based ROMs is independent of the quality factor.

\begin{figure}[!ht]
\centering
\includegraphics[width=0.45\textwidth]{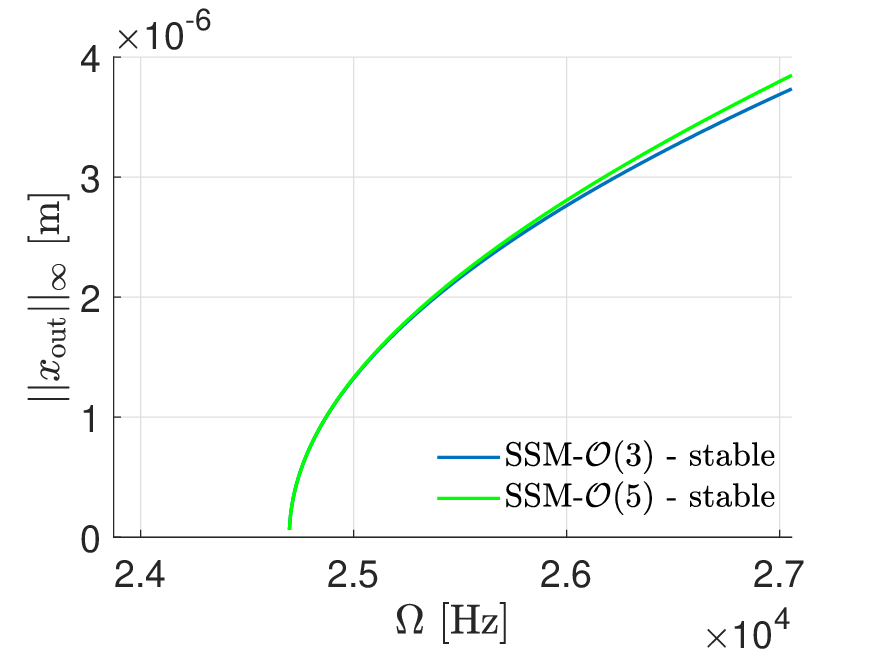}
\caption{Backbone curve for the first vibration mode of the MEMS device. The curve is obtained using SSM-based reduction truncated at cubic $\mathcal{O}(3)$ and fifth $\mathcal{O}(5)$ orders.}
\label{fig:mems-bb}
\end{figure}

We take the first mode as the master subspace and perform model reduction on the associated SSM. The first natural frequency is given as $\omega_1\approx24698$ Hz. We perform the non-intrusive computation of the autonomous SSM up to $\mathcal{O}(5)$. The associated backbone curve in terms of the $x_\mathrm{out}$ in Fig.~\ref{fig:memsplot} is shown in Fig.~\ref{fig:mems-bb}, from which we see the backbone curve converges well at $\mathcal{O}(3)$ approximation for $||x_\mathrm{out}||_\infty\leq 3~[\mu\mathrm{m}]$. In addition, we infer from Fig.~\ref{fig:mems-bb} that the nonlinear vibration around the first mode displays hardening behavior.

\begin{figure}[!ht]
\centering
\includegraphics[width=0.45\textwidth]{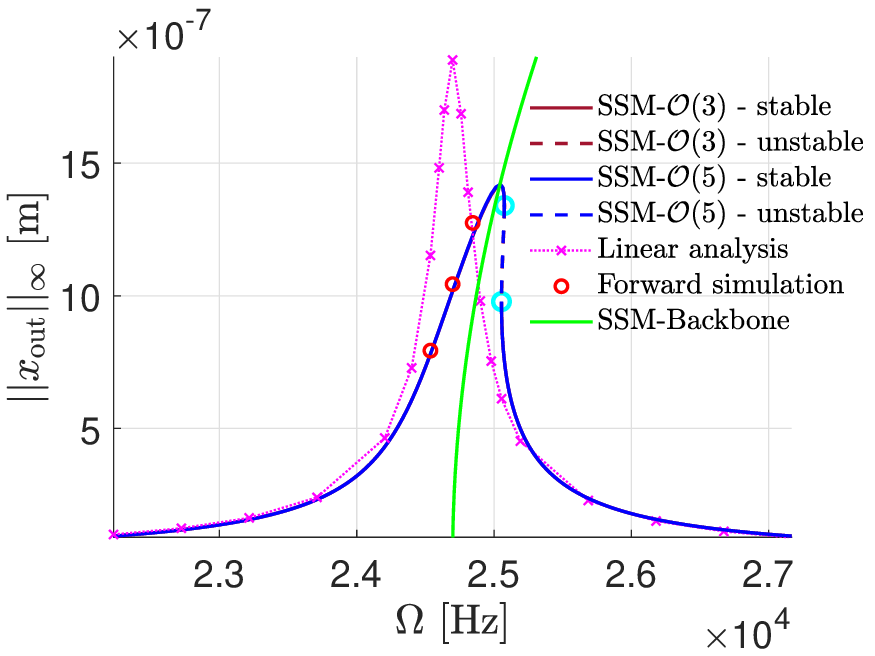}
\includegraphics[width=0.45\textwidth]{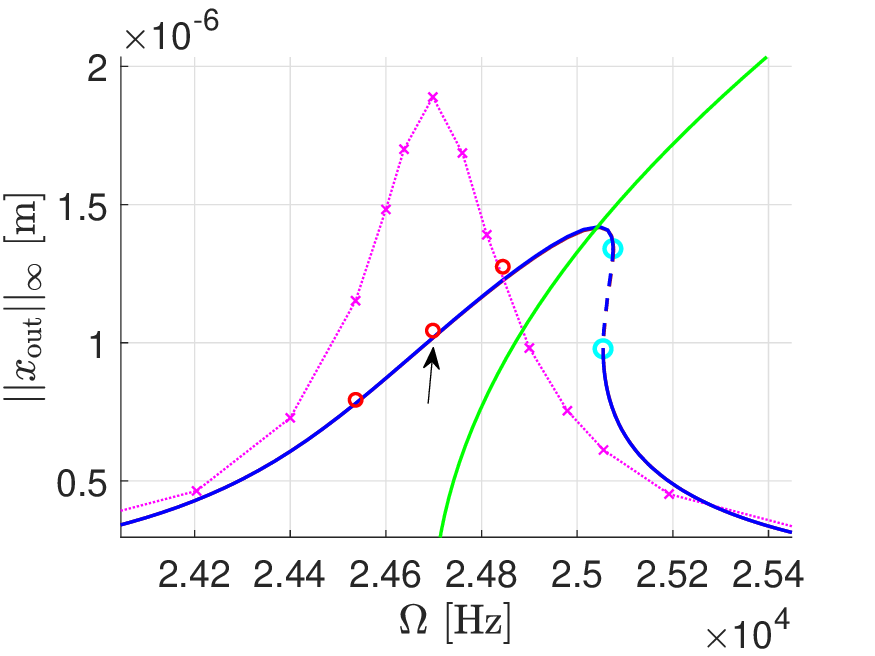}
\caption{Forced response curve for the primary resonance of the first vibration mode of the MEMS device. The right panel is a zoomed plot of the left panel around the resonance. Here, the results from cubic order $\mathcal{O}(3)$ SSM computations are very close to that of the fifth order $\mathcal{O}(5)$ computations and hence, hence, indistinguishable. The dotted line with cross markers denotes the forced responses from a linear analysis applied to the full system. Red circles represent results from forward simulations of the full system. The backbone curve at $\mathcal{O}(5)$ truncation at Fig.~\ref{fig:mems-bb} is also plotted here.}
\label{fig:mems-frc}
\end{figure}

Next, we add a concentrated harmonic forcing $\epsilon\cos\Omega t$ with $\epsilon=3.3386~\mu\mathrm{N}$ along the $x_\mathrm{out}$ direction (cf.~Fig.~\ref{fig:mems-frc}) and extract the forced response curve of the system around the first vibration mode. We perform a similar analysis as the one in Fig.~\ref{fig:shell_FRC} and find that it is sufficient to use TI-SSM solutions. So, we will only report TI-SSM solutions in the rest of this example. The FRC obtained using $\mathcal{O}(3)$ and $\mathcal{O}(5)$ is presented in Fig.~\ref{fig:mems-frc}. We note that the FRC predicted from the $\mathcal{O}(5)$ truncation matches excellently with that of the $\mathcal{O}(3)$ truncation such that the FRC at $\mathcal{O}(3)$ is indistinguishable. Along the FRC, two saddle-node (SN) bifurcated periodic orbits are detected, and the segment of periodic orbits is unstable between the two SN points. We also plot the FRC from a linear analysis applied to the full system in Fig.~\ref{fig:mems-frc}. As expected, the results from the linear analysis match well with that of the SSM-based nonlinear predictions when $\Omega$ is away from the resonance frequency, and a significant difference is observed for $\Omega\approx\omega_1$. In particular, the linear analysis overestimates the peak amplitude of the FRC and makes an incorrect prediction of the resonance frequency corresponding to the peak. We also plot the backbone curve (green line) along with the FRC in Fig.~\ref{fig:mems-frc}, from which we see that the peak on the FRC stays closely on the backbone curve.

\begin{figure}[!ht]
\centering
\includegraphics[width=0.45\textwidth]{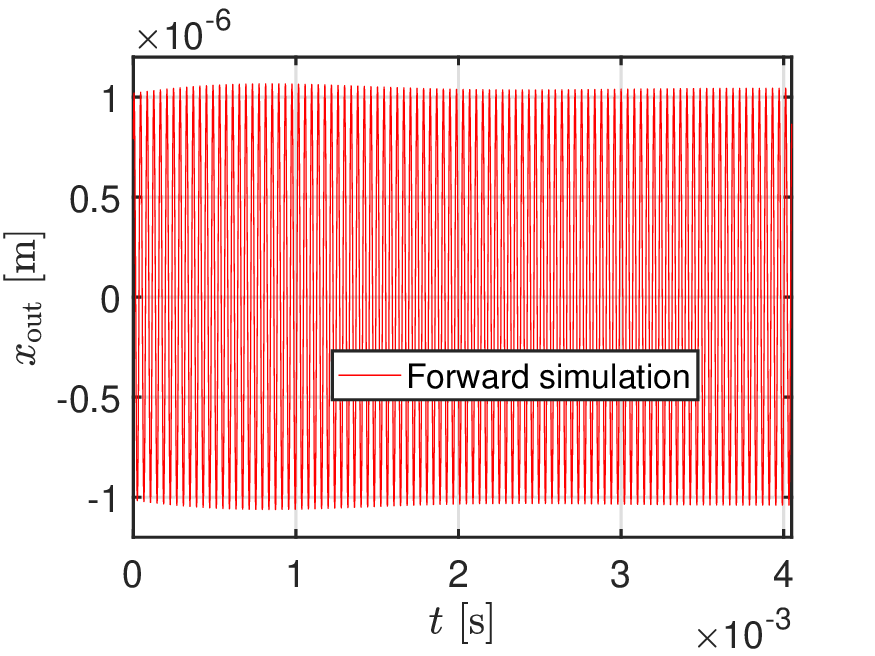}
\includegraphics[width=0.45\textwidth]{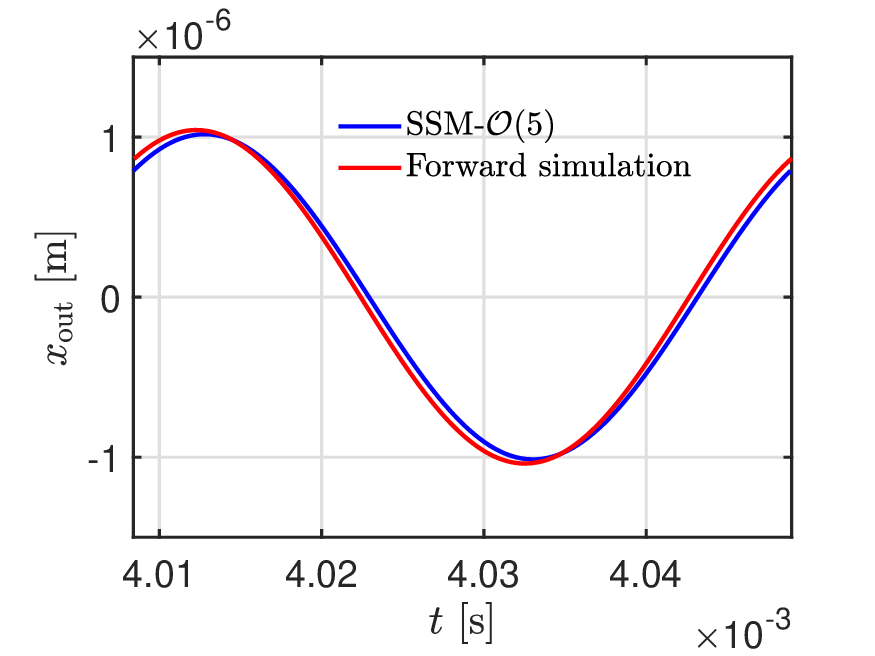}
\caption{Validation of the forced periodic orbit of the MEMS device using forward simulations of the full system. Here, the periodic orbit corresponds to the second circle in Fig.~\ref{fig:mems-frc}. In the left panel, the time history of the forward simulation in 100 cycles is presented. This simulation is initialized with a starting point on our predicted periodic orbit. In the right panel, we plot the response of forward simulation in the last circle along with our SSM-based prediction.}
\label{fig:mems-forward-traj}
\end{figure}

We now validate the effectiveness of the SSM-based predictions using forward simulations of the full system in COMSOL Multiphysics. Specifically, we take three sampled periodic orbits in the nonlinear regime of the FRC (cf.~the three red circles in Fig.~\ref{fig:mems-frc}) to perform the validation. We initialize each simulation with a starting point on the periodic orbit predicted using the SSM-based reduction and perform the simulation with 100 forcing cycles. Indeed, the simulation time of 100 cycles is long enough such that the simulation reaches a steady state, as illustrated in the left panel of Fig.~\ref{fig:mems-forward-traj}. In addition, we observe from the right panel of Fig.~\ref{fig:mems-forward-traj} that the periodic response in a steady state remains close to the SSM-based prediction. We take the magnitude of the signal in the steady state as the amplitude of the periodic orbit of the forward simulation. As seen in Fig.~\ref{fig:mems-frc}, the amplitude of three sampled periodic orbits matches well with that of our predicted periodic responses, demonstrating the SSM-based predictions' effectiveness.

In this example, the estimated memory cost for the non-intrusive algorithm is around 5 Gb, independently of the expansion orders of the SSM. The computational time at $\mathcal{O}(3)$ and $\mathcal{O}(5)$ computations are 2 hours and 19 hours. The averaged simulation for the three simulations is 89 hours. Again, we obtain a significant speed-up gain. We note the forward simulations are initialized with points on the SSM-predicted periodic orbits. In other words, the initial states are very close to the periodic orbits in a steady state. This significantly reduces the computational times of the forward simulations. Indeed, a much higher speed-up gain will be obtained if these forward simulations are initialized with states that are away from the periodic orbits in the steady state.

\section{Conclusion}
\label{sec:conclusion}

We have proposed a non-intrusive algorithm based on the stiffness evaluation procedure (STEP) ~\cite{Muravyov2003DeterminationStructures} for the automated computation of arbitrary dimensional spectral submanifolds (SSMs) up to arbitrary orders of accuracy. The algorithm works for general mechanical systems with up to cubic order nonlinearities, possibly with asymmetric damping and stiffness matrices, velocity-dependent nonlinear internal forces, and internal resonances. Our results enable rigorous, dynamics-based model reduction of nonlinear finite element (FE) models constructed via generic FE software. We have illustrated the effectiveness and efficiency of the non-intrusive algorithm towards nonlinear steady-state computation and bifurcation analysis over several examples of varying complexity, including a three-dimensional MEMS FE model with more than a million degrees of freedom.

We have implemented the non-intrusive algorithm in an open-source package, SSMTool~\cite{SSMTool2}, and importantly, established a coupling interface between SSMTool and a commercial finite element (FE) software, COMSOL Multiphysics. Such implementation significantly enhances the applicability of SSM-based model reduction for engineering applications, as illustrated in the last two examples considered in this study. We believe the implementation documented in this study paves the way for developing coupling with other FE codes and effectively bridges the gap between the rigorous theory of SSM-based exact model reduction and practical applications.

We conclude this study with a few ongoing and future research directions. Extension of the non-intrusive algorithm established here for constrained mechanical systems whose equations of motion are in the form of differential-algebraic equations (DAEs) is straightforward. We have used leading-order approximation for the non-autonomous part of SSMs throughout this study. Such an approximation is insufficient for systems with parametric excitation, where one has to compute higher-order expansions of the non-autonomous part of SSMs. The results of these ongoing projects will be published elsewhere.

\appendix
\section{Additional example: a cantilevered pipe conveying fluid}
\label{sec:add-example}

We consider a geometrically nonlinear viscoelastic cantilevered pipe conveying fluid subject to a harmonic base excitation $\epsilon\cos\Omega t$~\cite{li2023nonlinear}. Due to the flow-induced gyroscopic and follower forces, both the damping and stiffness matrices of the system are asymmetric. In addition, the system admits nonlinear damping up to cubic order. Such a nonlinear damping results from the viscoelastic constitutive law and large deformation. More details about the equations of motion of the pipe system can be found in Sect. 2 of~\cite{li2023nonlinear}.

SSM-based model reduction has been performed for the pipe system using the intrusive approach in~\cite{li2023nonlinear}. The obtained SSM-based ROMs predict the pipe system's free and forced vibrations accurately and efficiently. Remarkably, the ROMs enable predictions of global bifurcations such as saddle-node bifurcation on invariant circles and homoclinic bifurcations~\cite{li2023nonlinear}. 

Here, we will apply the proposed non-intrusive algorithms to the pipe system to demonstrate their effectiveness for generic systems with asymmetric damping, stiffness matrices, and velocity-dependent nonlinearities. In particular, we will compute the forced response of the pipe system in the post-flutter regime. We take the dimensionless flow velocity to be 6 and extract the forced response curve (FRC) for $\epsilon=0.0043$.

\begin{figure}[!ht]
\centering
\includegraphics[width=.45\textwidth]{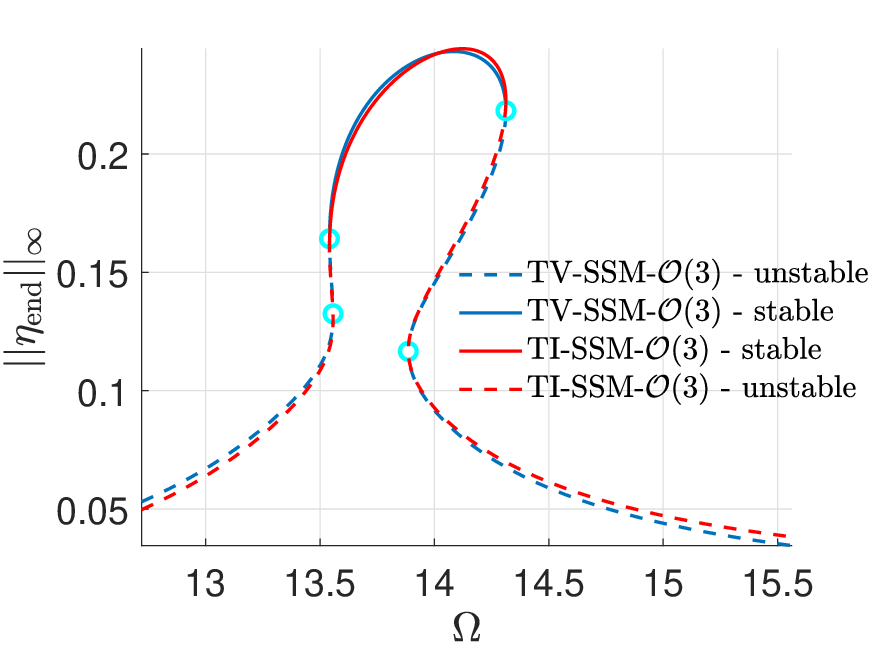}
\includegraphics[width=.45\textwidth]{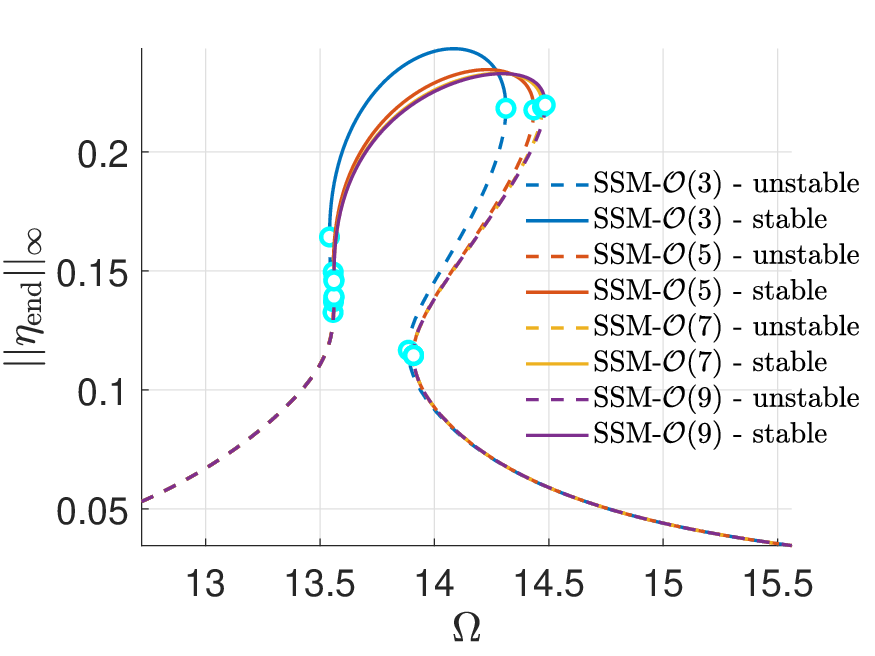}\\
\includegraphics[width=.45\textwidth]{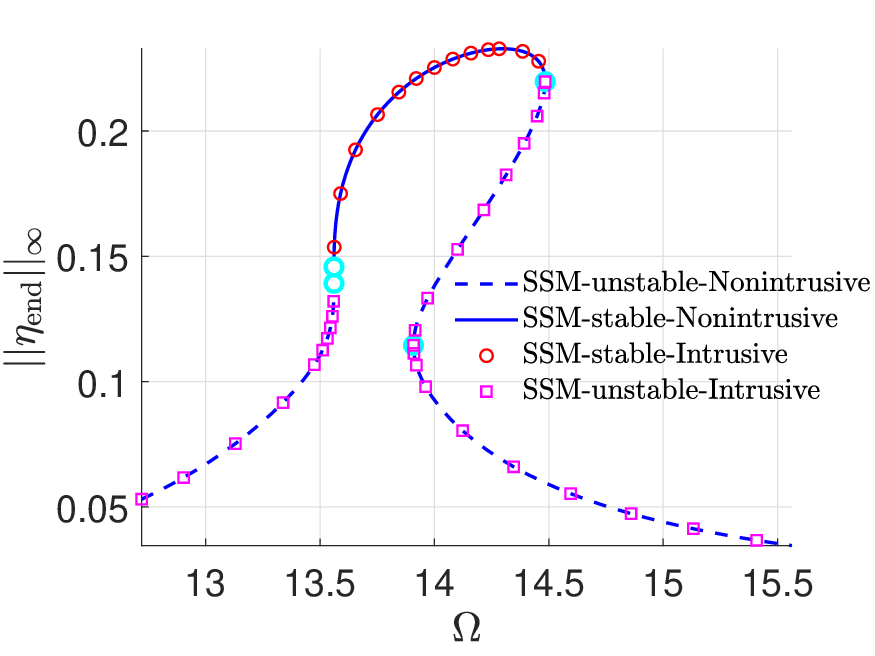}
\caption{Forced response curve regarding the tip deflection of the cantilevered pipe system. In the upper-left panel, we present the results obtained from $\mathcal{O}(3)$ approximation with and without the contribution of the leading-order non-autonomous part of SSM. Specifically, TV (TI) stands for time-varying (time-independent) solution with (without) the contribution. The results predicted from SSM-based reductions truncated at various expansion orders are shown in the upper-right panel. The lower panel presents the results at $\mathcal{O}(9)$ expansion, obtained using both the intrusive computation in~\cite{li2023nonlinear} and our proposed non-intrusive algorithm.}
\label{fig:pipe_frcs}
\end{figure}

In all the previous examples in Sect.~\ref{sec:examples}, we have safely ignored the contributions of the leading-order non-autonomous part when we extract the forced response curves (cf.~Remark~\ref{rk:ti-vs-tv}). However, as shown in the upper-left panel of Fig.~\ref{fig:pipe_frcs}, the difference between the TI-SSM and TV-SSM solutions is non-negligible and we have to account for the contribution of the non-autonomous part. Therefore, we use TV-SSM solutions for the rest of this example.

Next, we determine the expansion order based on the convergence of the forced response curve. As seen in the upper-right panel of Fig.~\ref{fig:pipe_frcs}, the FRC predicted from $\mathcal{O}(7)$ coincides well with that of the $\mathcal{O}(9)$. Thus, the forced response is converged at $\mathcal{O}(9)$, and we use the expansion at $\mathcal{O}(9)$ approximation, just as we did in~\cite{li2023nonlinear}. We plot the converged FRC at the lower panel of Fig.~\ref{fig:pipe_frcs}, along with the results computed via an intrusive computation and reported in the fourth panel of Fig.~8 in~\cite{li2023nonlinear}. We observe that the two schemes yield the same results, demonstrating that the non-intrusive approaches work for generic dynamical systems with asymmetric damping and stiffness matrices.

\section*{Acknowledgements}
ML acknowledges the financial support from the National Natural Science Foundation of China (No. 12302014) and Shenzhen Science and Technology Innovation Commission (No. 20231115172355001).

\appendix

\bibliography{references} 
\bibliographystyle{elsarticle-num} 
\end{document}